\documentclass[10pt,a4paper]{article}
\usepackage[T1]{fontenc}
\usepackage[utf8]{inputenc}
\usepackage{lmodern}
\usepackage{amsmath,amssymb,amsthm,mathtools}
\usepackage[margin=27mm,headheight=14pt]{geometry}
\usepackage{microtype}
\usepackage{needspace}
\usepackage[dvipsnames]{xcolor}
\usepackage{hyperref}
\hypersetup{colorlinks=true,linkcolor=MidnightBlue,citecolor=MidnightBlue,urlcolor=MidnightBlue,pdftitle={Radial Symmetry and Nondegeneracy for Singular Liouville Systems with Small Singular Exponents}}
\numberwithin{equation}{section}
\newtheorem{theorem}{Theorem}[section]
\newtheorem{proposition}[theorem]{Proposition}
\newtheorem{lemma}[theorem]{Lemma}
\newtheorem{corollary}[theorem]{Corollary}
\newtheorem{knownthm}{Theorem}

\DeclareMathOperator{\Ker}{Ker}
\DeclareMathOperator{\Ran}{Ran}

\theoremstyle{definition}
\newtheorem{definition}[theorem]{Definition}
\theoremstyle{remark}
\newtheorem{remark}[theorem]{Remark}
\newtheorem{example}[theorem]{Example}

\newcommand{\R}{\mathbb R}
\newcommand{\dd}{\,\mathrm d}
\newcommand{\e}{\mathrm e}
\DeclareMathOperator{\dist}{dist}
\DeclareMathOperator{\Span}{span}
\DeclareMathOperator{\supp}{supp}
\title{Radial Symmetry and Nondegeneracy for Singular Liouville Systems with Small Singular Exponents}
\author{
    Zetao Cheng
    \and
	Haoyu Li
	\and
	Lei Zhang
}

\date{}

\begin{document}
\maketitle

\begin{abstract}
We extend the classification (Lin and Zhang, Ann. Inst. H. Poincar\'e Anal. Non Lin\'eaire \textbf{27}, 2010) to finite-mass singular Liouville systems with sufficiently small positive singular exponents. For symmetric, entrywise nonnegative, irreducible and invertible coupling matrices, every solution whose mass vector lies near a fixed compact subset of the regular mass surface is radial about the singular point and is uniquely determined by its component masses up to dilation. We also establish nondegeneracy modulo dilation. To overcome the failure of the moving-plane argument for positive exponents, we use compactness and integral identities for bounded solutions of the linearized system. The results also hold for nonzero nonnegative exponent vectors and for small mixed-sign vectors satisfying an explicit integral inequality.
\end{abstract}

\tableofcontents

\section{Introduction}\label{sec:introduction}

We consider the singular Liouville system
\begin{equation}\label{eq:system}
-\Delta u_i=\sum_{j\in I}a_{ij}|x|^{\beta_j}\e^{u_j}
\quad\text{in }\R^2,\qquad i\in I,
\end{equation}
subject to the finite-mass conditions
\begin{equation}\label{eq:mass}
\sigma_i(u)=\frac{1}{2\pi}\int_{\R^2}|x|^{\beta_i}\e^{u_i(x)}\dd x<\infty,
\qquad i\in I.
\end{equation}
Here $I=\{1,\ldots,N\}$, $N\geq2$ and $\beta=(\beta_1,\ldots,\beta_N)\in(-2,\infty)^N$ is the singular exponent vector. The real coupling matrix $A=(a_{ij})_{i,j\in I}$ satisfies
\begin{equation*}\tag{H1}\label{eq:H1}
A \text{ is symmetric, entrywise nonnegative, irreducible and invertible.}
\end{equation*}
Irreducibility means that no nonempty proper subset $J\subset I$ satisfies $a_{ij}=0$ for all $i\in J$ and $j\notin J$. We do not assume that $A$ is positive definite.

A solution is a vector in $C^0(\R^2,\R^N)\cap C^2(\R^2\setminus\{0\},\R^N)$ satisfying \eqref{eq:system} in the sense of distributions. Thus $u_i$ is bounded at the origin. The prescribed singularity is contained in the weight $|x|^{\beta_i}$.

For the scalar Liouville equation, Chen and Li \cite{CL1991}
classified finite-mass entire solutions.
Prajapat and Tarantello \cite{PT2001} obtained classification
results for the corresponding equation with a radial power weight.
For singular sources of integer multiplicity, del Pino, Esposito
and Musso \cite{DEM2012} proved that all bounded solutions of
the linearized equation are generated by differentiation of
the explicit solution family.

For the regular system, Chipot, Shafrir and Wolansky proved radial symmetry about a common center \cite{CSW1997}. Lin and Zhang proved fixed-mass uniqueness \cite{LZ2010} and identified the bounded kernel at a regular radial solution as dilation and the two common translations \cite[Theorem~2.1(2)]{LZ2013}. For singular exponents, they classified radial solutions and proved fixed-mass uniqueness \cite[Theorem~1.1]{LZ2014}. They also proved radial symmetry when all exponents are nonpositive and not all zero and nondegeneracy when all exponents are strictly negative \cite[Theorems~1.2--1.3]{LZ2014}.

For the $SU(n+1)$ Toda system, Lin, Wei and Ye \cite{LWY2012}
classified finite-mass solutions with a singular source and proved
nondegeneracy of the linearized equations.

Li--Zhang \cite{LiZhang2003} developed the method of moving
spheres to prove Liouville theorems and Harnack inequalities for
semilinear elliptic equations with critical growth.
Wei and Xu \cite{WX1999} obtained classification results for
higher order conformally invariant equations.
In their study of the prescribing $Q$-curvature problem on $S^n$,
they proved that every bounded solution of the linearized integral
equation at a standard solution is a linear combination of the
dilation and translation derivatives \cite{WX2009}.
For singular fourth-order equations, Ahmedou, Wu and Zhang
\cite{AWZ2022} obtained classification and quantization results
and proved a priori estimates for the prescribing $Q$-curvature
problem on compact four-manifolds.
Rigidity can also be obtained from integral identities of Obata type.
Lin and Ma \cite{LM2024} proved that positive solutions of a class
of quasilinear $p$-Laplace type equations on the sphere are constant
in a subcritical range.
For a survey of Obata type integral identities and their applications,
we refer to Lin, Ma and Ou \cite{LMO2024}.

Busca and Sirakov \cite{BS2000} proved radial symmetry for
positive solutions of cooperative elliptic systems that tend to
zero at infinity under suitable structural assumptions.
Rigidity results for elliptic systems also concern equality of components.
Quittner and Souplet \cite{QS2012} gave sufficient conditions under which
positive entire solutions of a semilinear elliptic system have equal
components and obtained classification and Liouville theorems.
For systems with power nonlinearities, Jing et al.
\cite{JLLLW2024} obtained uniqueness and multiplicity results
for positive solutions of $n$-coupled elliptic systems.
For equal-frequency two-component cubic Schr\"odinger systems
in dimensions two and three, Chen, Liu, Wei and Yang
\cite{CLWY2026} proved uniqueness of positive finite-energy
solutions up to simultaneous translations in the weak positive coupling range.

The connection between Liouville theorems and a priori estimates for
elliptic equations and systems was developed by
Pol\'a\v{c}ik, Quittner and Souplet \cite{PQS2007}.
The exponential integrability and local estimates of Brezis
and Merle \cite{BM1991} provide basic tools for the blow-up
analysis of two-dimensional equations with exponential nonlinearities.

For scalar equations associated with weighted interpolation inequalities,
Dolbeault, Esteban and Loss \cite{DEL2016} obtained rigidity results
and determined the symmetry range in the Caffarelli--Kohn--Nirenberg
inequalities they considered.
Dolbeault, Muratori and Nazaret \cite{DMN2017} proved radial symmetry
of optimizers when the exponent of the power weight is sufficiently
small, using compactness and a perturbation argument based on a
spectral gap.

We consider small exponent vectors whose components are either all nonnegative or have different signs. Our aim is to prove radial symmetry about the prescribed singular point and identify the full bounded kernel, uniformly for masses near a compact subset of the regular mass surface.

Set
\begin{equation}\label{eq:gamma-definitions}
\gamma_i=2+\beta_i,\qquad m=A\sigma,
\end{equation}
and, for every nonempty $J\subset I$, define
\begin{equation}\label{eq:lambda}
\Lambda_J^\beta(\sigma)
=2\sum_{i\in J}(2+\beta_i)\sigma_i
-\sum_{i,j\in J}a_{ij}\sigma_i\sigma_j.
\end{equation}
The admissible radial mass surface is
\begin{equation}\label{eq:mass-surface}
\Gamma_\beta=\left\{\sigma\in(0,\infty)^N:
\Lambda_I^\beta(\sigma)=0,\quad
\Lambda_J^\beta(\sigma)>0\ \text{for }\varnothing\neq J\subsetneq I
\right\}.
\end{equation}
We write $\Gamma_0$ for the regular surface. Distances between mass vectors are Euclidean.

The dilation
\begin{equation}\label{eq:dilation}
u_i^\lambda(x)=u_i(\lambda x)+(2+\beta_i)\log\lambda,
\qquad \lambda>0,
\end{equation}
preserves \eqref{eq:system} and every component mass. We fix this freedom by the total-mass median condition
\begin{equation}\label{eq:median}
\sum_i\int_{B_1}|x|^{\beta_i}\e^{u_i(x)}\dd x
=\frac12\sum_i\int_{\R^2}|x|^{\beta_i}\e^{u_i(x)}\dd x.
\end{equation}
The radial classification theorem of Lin and Zhang \cite[Theorem~1.1]{LZ2014} implies that every $\eta\in\Gamma_0$ determines a centered regular radial solution uniquely up to dilation. We denote by $U^\eta$ the representative satisfying \eqref{eq:median}.

By a bounded solution of the linearized system
\begin{equation}\label{eq:linearized}
-\Delta\phi_i=\sum_j a_{ij}|x|^{\beta_j}\e^{u_j}\phi_j
\end{equation}
we mean $\phi\in L^\infty(\R^2,\R^N)$ satisfying the equation in the sense of distributions on the whole plane. For any solution $u$, dilation generates the solution $Z$ of \eqref{eq:linearized}, given by
\begin{equation}\label{eq:dilation-field}
Z_i(x)=x\cdot\nabla u_i(x)+2+\beta_i.
\end{equation}

Our first theorem treats a nonzero exponent vector with small nonnegative components. The solution's mass vector is required to lie near a fixed compact subset of $\Gamma_0$.

\begin{theorem}\label{thm:nonnegative}
Assume \eqref{eq:H1} and let $K\subset\Gamma_0$ be nonempty and compact.
There exist $\varepsilon,\delta_0>0$ such that every solution of
\eqref{eq:system}--\eqref{eq:mass} with
\[
 \beta\in[0,\varepsilon)^N\setminus\{0\},\qquad
 \dist(\sigma(u),K)<\delta_0
\]
satisfies
\[
 u_i=u_i(|x|),\qquad \sigma(u)\in\Gamma_\beta,\qquad
 \{\phi\in L^\infty:\phi\text{ solves }\eqref{eq:linearized}\}
       =\Span\{Z\}.
\]
\end{theorem}

To allow exponent vectors with components of different signs, we impose the following
condition on the integrals $\int_{B_r}\e^{U_i^\eta(x)}\dd x$.

For $\eta\in\Gamma_0$, let
\begin{equation}\label{eq:disk-integral}
E_i^\eta(r)=\int_{B_r}\e^{U_i^\eta(x)}\dd x,
\qquad r>0.
\end{equation}
Here $U^\eta$ solves \eqref{eq:system} with $\beta=0$ and $(2\pi)^{-1}\int_{\R^2}\e^{U_i^\eta}=\eta_i$ for $i\in I$. Let $\mathcal B$ be a nonempty compact subset of the unit $\ell^1$ sphere,
\[
\sum_i|b_i|=1\qquad (b\in\mathcal B).
\]
We assume that there exists a constant $\kappa>0$, depending only on $A$, $K$ and $\mathcal B$, such that
\begin{equation*}\tag{NC}\label{eq:NC}
\left|\sum_i b_iE_i^\eta(r)\right|
\geq\kappa\sum_i|b_i|E_i^\eta(r)
\quad\text{for every }\eta\in K,\ b\in\mathcal B,\ r>0.
\end{equation*}
This condition is invariant under a common dilation of the centered profile. It holds with $\kappa=1$ when every component of $b$ is nonnegative. Section~\ref{sec:radial} constructs compact sets of vectors $b$ with both positive and negative components for which it also holds.

\begin{theorem}\label{thm:signed}
Assume \eqref{eq:H1} and let $K\subset\Gamma_0$ and $\mathcal B$ be
nonempty compact sets satisfying \eqref{eq:NC}.
There exist $\varepsilon,\delta_0>0$ such that every solution of
\eqref{eq:system}--\eqref{eq:mass} with
\[
 \beta=\tau b,\qquad 0<\tau<\varepsilon,\qquad b\in\mathcal B,
 \qquad \dist(\sigma(u),K)<\delta_0
\]
satisfies
\[
 u_i=u_i(|x|),\qquad \sigma(u)\in\Gamma_\beta,\qquad
 \{\phi\in L^\infty:\phi\text{ solves }\eqref{eq:linearized}\}
       =\Span\{Z\}.
\]
\end{theorem}

By \cite[Theorem~1.1]{LZ2014}, every mass in $\Gamma_\beta$ is realized by a centered radial solution, unique up to \eqref{eq:dilation}. Thus, in either parameter range of Theorems~\ref{thm:nonnegative} and~\ref{thm:signed}, solutions with the same exponents and component masses differ only by dilation.

Because $\sum_i|b_i|=1$, every $b_i\geq-1$. We may decrease $\varepsilon$ so that $\varepsilon<2$. Hence all exponents in Theorem~\ref{thm:signed} satisfy $\beta_i>-2$.

The proof starts with compactness for solutions normalized by
\eqref{eq:median}. Under \eqref{eq:NC}, Lemma~\ref{lem:rectangle}
and the integral identities exclude concentration and recover all
component masses. The regular classification initially determines
the limit only up to translation and dilation. We locate its center
by dividing \eqref{mom:identity-eq} by $\tau$ before passing to the
limit. The resulting identity is combined with the exact formula
\begin{equation}\label{eq:center-preview}
\int_{\R^2}\frac{x}{|x|^2}\sum_i b_i\e^{U_i^\eta(x-a)}\dd x
=\frac{a}{|a|^2}\sum_i b_iE_i^\eta(|a|),
\qquad a\neq0,
\end{equation}
whose right-hand side is nonzero by \eqref{eq:NC}. The same formula
applies after a common dilation, so the limit is centered at zero.
Condition \eqref{eq:median} then determines its dilation.

We next determine the bounded solutions of the linearized system at
an arbitrary solution under consideration. In a contradiction
argument, we subtract a multiple of $Z$ and normalize the supremum
norm. The regular kernel theorem restricts the limit to the two
translation solutions. The integral identity in
Lemma~\ref{rig:linearized-moment} excludes these terms, since
\eqref{eq:NC} as $r\downarrow0$ bounds
$\left|\sum_i b_i\e^{U_i^\eta(0)}\right|$ away from zero.
A scalar comparison after inversion then gives uniform convergence
to zero on the plane, contradicting the normalization.
Thus the bounded kernel is spanned by $Z$. The rotation derivative
belongs to this kernel and vanishes at the origin, which proves
radial symmetry. The published radial classification gives the mass
condition and uniqueness up to dilation.

The study of radial symmetry for small positive singularities is
motivated by the scalar $Q$-curvature work of Li, Li, Wei and Yang
\cite[Sections~6--7]{QDraft}.

Section~\ref{sec:published} recalls the classical results and proves
smooth dependence of the radial solutions on their initial values
and component masses. It fixes dilation by \eqref{eq:median} and
obtains uniform estimates used to construct vectors satisfying
\eqref{eq:NC}. Section~\ref{sec:moments} proves the integral identity
for the original system and the formula \eqref{eq:center-preview}.
Section~\ref{sec:compactness} proves compactness for normalized
sequences, preserves all component masses and identifies the center
of the regular limit. Section~\ref{sec:rigidity} establishes uniform
estimates after inversion, proves the integral identity for bounded
solutions of the linearized system and determines their kernel.
It then proves the main theorems by applying this conclusion to the
rotation derivative.

\section{Refined analysis on the radial solutions}\label{sec:published}\label{sec:radial}

\subsection{Preliminary results}

We first recall the classification, asymptotic estimates, and kernel results used below, and then construct the normalized radial families and their mass coordinates. Throughout this section, $A$ satisfies \eqref{eq:H1}.

\begin{knownthm}\label{thm:published-classification}
Assume \eqref{eq:H1}. For fixed $\beta_i>-2$, a positive mass vector $\sigma$ is realized by a centered radial solution of \eqref{eq:system} if and only if $\sigma\in\Gamma_\beta$. Two centered radial solutions with the same component masses differ by a dilation. When $\beta=0$, every finite-mass entire solution is radial about a common center, and every component decreases strictly away from that center.
\end{knownthm}

\noindent
Lin and Zhang proved radial realization and fixed-mass uniqueness in \cite[Theorem~1.1]{LZ2014}. The regular common-center assertion is the theorem of Chipot, Shafrir, and Wolansky \cite{CSW1997}, quoted as Theorem~A in \cite{LZ2010}. Regular fixed-mass uniqueness is also contained in \cite[Theorem~1.1]{LZ2010}. These results use precisely the matrix assumptions in \eqref{eq:H1}.

\begin{knownthm}\label{thm:published-asymptotics}
Let $u$ be a solution of \eqref{eq:system}--\eqref{eq:mass}, and set $m=A\sigma$. Then
\[
m_i>2+\beta_i\qquad(i\in I).
\]
There exist constants $c_i$ and $\delta_u>0$ such that as $|x|\to\infty$,
\begin{equation}\label{eq:published-asymptotics}
u_i(x)=-m_i\log|x|+c_i+O(|x|^{-\delta_u}),
\qquad
\nabla u_i(x)=-m_i\frac{x}{|x|^2}+O(|x|^{-1-\delta_u}).
\end{equation}
Moreover, for suitable constants $C_i$,
\begin{equation}\label{eq:published-potential}
u_i(x)=C_i+\frac1{2\pi}\sum_j a_{ij}
\int_{\R^2}\bigl[-\log|x-y|+\log(1+|y|)\bigr]
|y|^{\beta_j}\e^{u_j(y)}\dd y.
\end{equation}
The constants $c_i$, $C_i$, and $\delta_u$ may depend on the solution $u$.
\end{knownthm}

\noindent
The value and gradient estimates follow from \cite[Lemma~2.1]{LZ2014}. The integral representation is proved there in Section~6, equations (42)--(43) of the author manuscript linked in the bibliography.

\Needspace{10\baselineskip}
The following kernel descriptions are due to Lin and Zhang
\cite[Theorem~2.1(2)]{LZ2013} and
\cite[Lemma~3.2 and Theorem~1.3]{LZ2014}.

\begin{knownthm}\label{thm:published-regular-kernel}
Assume \eqref{eq:H1}, and let $U$ be a finite-mass radial solution of \eqref{eq:system}.
\begin{enumerate}
\item[(i)] If $\beta=0$, every bounded entire solution of
\begin{equation}\label{eq:regular-linearized}
-\Delta\phi_i=\sum_j a_{ij}\e^{U_j}\phi_j
\end{equation}
has the form
\begin{equation}\label{eq:regular-kernel}
\phi_i(r,\theta)=c_0\bigl(rU_i'(r)+2\bigr)
+c_1U_i'(r)\cos\theta+c_2U_i'(r)\sin\theta,
\end{equation}
where $c_0,c_1,c_2\in\R$ do not depend on $i$.
\item[(ii)] If $\beta\in(-2,\infty)^N$, every bounded radial solution of \eqref{eq:linearized} at $U$ has the form
\begin{equation}\label{eq:singular-radial-kernel}
\phi_i(r)=c_0\bigl(rU_i'(r)+2+\beta_i\bigr),
\end{equation}
where $c_0\in\R$ does not depend on $i$.
If, in addition, $-2<\beta_i<0$ for every $i\in I$, every bounded entire solution of \eqref{eq:linearized} has the form \eqref{eq:singular-radial-kernel}.
\end{enumerate}
\end{knownthm}

\begin{remark}\label{rem:published-scope}
The asymptotic estimates in Theorem~\ref{thm:published-asymptotics} hold for each fixed solution, with constants that may depend on that solution. Uniform estimates for the radial parameter families are proved below. Those for normalized solution sequences are proved in Section~\ref{sec:rigidity}.
\end{remark}

We now consider the dependence of radial solutions on their initial values
and masses. The logarithmic radial variable will also be used in the
linearized equation.

For a radial solution $U$, we introduce \emph{the Emden-Fowler transform}
\begin{equation}\label{eq:radial-cylinder}
 t=\log r,
 \qquad w_i(t)=U_i(\e^t),
 \qquad q_i(t)=\e^{\gamma_i t+w_i(t)}.
\end{equation}
Then a routine computation gives
\begin{equation}\label{eq:radial-ode}
 w_i''=-\sum_j a_{ij}q_j,
 \qquad
 \sigma_i=\int_{\R}q_i(t)\dd t.
\end{equation}
The normalization of $\sigma_i$ already removes the polar factor $2\pi$. We also write $m=A\sigma$ and $d=m-\gamma$.

The following lemma will be used in the compactness argument to rule out concentration of only part of the mass and to show that every component retains its full mass in the limit. Recall that we select a compact set $K\subset\subset\Gamma_\beta$.

\begin{lemma}\label{lem:rectangle}
Let $\sigma\in\Gamma_\beta$ and $0\le\tau_i\le\sigma_i$. If $\tau\notin\{0,\sigma\}$, then $\Lambda_I^\beta(\tau)>0$. More precisely, with $p_i=\tau_i/\sigma_i$ and
\[
 \omega_J=\prod_{i\in J}p_i\prod_{i\notin J}(1-p_i),
\]
one has
\begin{equation}\label{eq:rectangle-identity}
 \Lambda_I^\beta(\tau)
 =\sum_{J\subseteq I}\omega_J\Lambda_J^\beta(\sigma)
  +\sum_i a_{ii}\tau_i(\sigma_i-\tau_i),
\end{equation}
where $\Lambda_\varnothing^\beta=0$. For $\beta=0$, if $K\subset\Gamma_0$ is compact, the positivity is uniform for $\sigma\in K$ whenever both $|\tau|$ and $|\sigma-\tau|$ are bounded below by a fixed positive number.
\end{lemma}

\begin{proof}
By the definition of $\omega_J$, we get
\[
\sum_{J\subseteq I}\omega_J
=\prod_{k\in I}\bigl(p_k+(1-p_k)\bigr)=1,
\]
and
\[
\sum_{\substack{J\subseteq I\\ i\in J}}\omega_J
=p_i\prod_{k\ne i}\bigl(p_k+(1-p_k)\bigr)
=p_i.
\]
Similarly, if $i\ne j$, then
\[
\sum_{\substack{J\subseteq I\\ i,j\in J}}\omega_J
=p_ip_j
\prod_{k\notin\{i,j\}}\bigl(p_k+(1-p_k)\bigr)
=p_ip_j.
\]
Together, they give
\[
\begin{aligned}
\sum_{J\subseteq I}\omega_J\Lambda_J^\beta(\sigma)
&=2\sum_i(2+\beta_i)\sigma_i
  \sum_{\substack{J\subseteq I\\ i\in J}}\omega_J-\sum_{i\ne j}a_{ij}\sigma_i\sigma_j
  \sum_{\substack{J\subseteq I\\ i,j\in J}}\omega_J-\sum_i a_{ii}\sigma_i^2
  \sum_{\substack{J\subseteq I\\ i\in J}}\omega_J\\
&=2\sum_i(2+\beta_i)\sigma_i p_i
  -\sum_{i\ne j}a_{ij}\sigma_i\sigma_jp_ip_j
  -\sum_i a_{ii}\sigma_i^2p_i\\
&=2\sum_i(2+\beta_i)\tau_i
  -\sum_{i,j}a_{ij}\tau_i\tau_j
  -\sum_i a_{ii}\tau_i(\sigma_i-\tau_i)\\
&=\Lambda_I^\beta(\tau)
  -\sum_i a_{ii}\tau_i(\sigma_i-\tau_i).
\end{aligned}
\]
Rearranging yields \eqref{eq:rectangle-identity}. Every term on its right is nonnegative. Since $p$ is neither $0$ nor $(1,\ldots,1)$, some nonempty proper $J$ has $\omega_J>0$, and its contribution is strictly positive. For $\sigma\in K$, let $c=\min_{\sigma\in K,\ \varnothing\ne J\subsetneq I}\Lambda_J^0(\sigma)>0$. By the the compactness of $K$ we find a $c$,
\[
\Lambda_J^0(\sigma)\ge c
\qquad
\text{for every }\sigma\in K
\text{ and }\varnothing\ne J\subsetneq I.
\]
Thus \eqref{eq:rectangle-identity} gives
\[
\begin{aligned}
\Lambda_I^0(\tau)&=\sum_{\varnothing\ne J\subsetneq I}
  \omega_J\Lambda_J^0(\sigma)
  +\sum_i a_{ii}\tau_i(\sigma_i-\tau_i)\ge
  \sum_{\varnothing\ne J\subsetneq I}
  \omega_J\Lambda_J^0(\sigma)\ge c\sum_{\varnothing\ne J\subsetneq I}\omega_J\\
  &\geq c \left(1-\prod_i p_i-\prod_i(1-p_i)\right).
\end{aligned}
\]
\end{proof}

\subsection{Dependence on the initial values}

To use the masses as local coordinates, we need smooth dependence of
radial solutions on the initial values, including estimates at both
ends of the cylinder.

Fix smooth functions $\chi_-+\chi_+=1$, with $\chi_+=0$ on
$(-\infty,0]$ and $\chi_+=1$ on $[1,\infty)$.
\begin{definition}
For an integer $\ell\ge0$ and $\delta>0$, we denote 
\[
 \|f\|_{C^\ell_\delta}
 =\max_{0\le j\le\ell}\sup_{t\in\R}
       \e^{\delta|t|}|\partial_t^j f(t)|.
\]
The space $C_\delta^\ell$ is the closure of $C_0^\infty$ in $\|\cdot\|_{C_\delta^\ell}$.
\end{definition}
In follows, we study the Emden--Fowler transform of $U=(U_1,\cdots,U_N)$. Here, $U$ is a radial solution solves
$-\Delta U_i=\sum_j a_{ij}|x|^{\beta_j}\e^{U_j}$ continuous at the origin and twice continuously differentiable away
from it. Finite component masses are required for the asymptotic estimates
as $r\to\infty$ and for the mass coordinates.

\begin{proposition}\label{prop:radial-parameters}
\begin{enumerate}
\renewcommand{\labelenumi}{(\roman{enumi})}

\item For every $\beta\in(-2,\infty)^N$ and $c\in\R^N$, there is a unique
radial solution with $U_i(0)=c_i$, defined at every finite radius.
It satisfies
\begin{align}
 w_i'(t)&=-\sum_j a_{ij}\int_{-\infty}^tq_j(s)\dd s,
 \label{eq:radial-first}\\
 w_i(t)&=c_i-\sum_j a_{ij}\int_{-\infty}^t(t-s)q_j(s)\dd s.
 \label{eq:radial-volterra}
\end{align}

\item Suppose that
\[
 \int_{\R}q_i(t,\beta_*,c_*)\,\dd t<\infty
 \qquad(i\in I).
\]
There exist $\rho,C,\nu>0$ such that
\begin{equation}\label{eq:radial-uniform-decay}
 0<q_i(t,\beta,c)\le C\e^{-\nu|t|}
 \qquad(t\in\R,\ i\in I),
\end{equation}
whenever $|\beta-\beta_*|+|c-c_*|<\rho$.

\item On this neighborhood, $q$ and all its parameter derivatives belong to
$C^\ell_\delta$ and depend continuously on $(\beta,c)$ in this norm,
for every fixed integer $\ell\ge0$ and every $0<\delta<\nu$.

\item In the same neighborhood, the masses, $m=A\sigma$, and
\begin{equation}\label{eq:right-constant}
 c^+=c+A\int_{\R}tq(t)\dd t
\end{equation}
depend smoothly on $(\beta,c)$.

\item For the same $\ell$ and $\delta$, the remainder
\begin{equation}\label{eq:weighted-smoothness}
 w-c\chi_--(c^+-mt)\chi_+
\end{equation}
and all its parameter derivatives belong to $C^\ell_\delta$
and depend continuously on $(\beta,c)$ in this norm.

\end{enumerate}
\end{proposition}

\begin{proof}
Write $p=(\beta,c)$.
\begin{enumerate}
\renewcommand{\labelenumi}{(\roman{enumi})}

\item \emph{Existence and local parameter dependence.}
The existence and uniqueness of the radial solution for every
$\beta_i>-2$ and $c\in\R^N$ follow from
\cite[Lemma~2.3]{LZ2014}. Its integral representation gives
\eqref{eq:radial-first} and \eqref{eq:radial-volterra}.
We also establish the local parameter dependence needed below.

Fix a finite-mass parameter pair $p_*=(\beta_*,c_*)$.
Throughout the remaining proof, the parameters vary in a sufficiently
small bounded neighborhood $\mathcal U$ of $p_*$, on which
\[
 \gamma_i=2+\beta_i\ge\gamma_*>0,
 \qquad |c_i|\le C_*.
\]
We make explicit the parameter dependence of the iteration used in
the proof of \cite[Lemma~2.3]{LZ2014}. In the $r$--variable, the
integral equation is
\[
 U_i(r)=c_i-\sum_j a_{ij}\int_0^r
       \log\frac rs\,s^{1+\beta_j}\e^{U_j(s)}\,\dd s.
\]
Write $h=U-c$ and define
\[
 (\mathcal T_p h)_i(r)
 =-\sum_j a_{ij}\int_0^r
 \log\frac rs\,s^{1+\beta_j}\e^{c_j+h_j(s)}\,\dd s,
 \qquad i\in I.
\]
For $\|h\|_\infty,\|\widetilde h\|_\infty\le1$, we have
\[
\begin{aligned}
 |(\mathcal T_p h)_i(r)-(\mathcal T_p\widetilde h)_i(r)|
 &\le
 \sum_j a_{ij}\int_0^r
 \log\frac rs\,s^{1+\beta_j}\e^{c_j}
 \bigl|\e^{h_j(s)}-\e^{\widetilde h_j(s)}\bigr|\,\dd s\\
 &\le
 \e^{C_*+1}\left(\sum_j a_{ij}\right)
 \|h-\widetilde h\|_\infty
 \int_0^r\log\frac rs\,s^{\gamma_*-1}\,\dd s\\
 &=
 \frac{\e^{C_*+1}}{\gamma_*^2}
 \left(\sum_j a_{ij}\right)
 r^{\gamma_*}\|h-\widetilde h\|_\infty.
\end{aligned}
\]
For sufficiently small $r_0\in(0,1)$, the operator preserves the
unit ball and is a contraction there, uniformly for $p\in\mathcal U$.

Parameter derivatives of the kernel introduce powers of $\log s$.
For every fixed integer $k\ge0$ and sufficiently small
$\varepsilon_0>0$,
\[
 \int_0^r
 \log\frac rs\,s^{\gamma_*-1}(1+|\log s|)^k\,\dd s
 \le C_{k,\varepsilon_0}r^{\gamma_*-\varepsilon_0}
 \qquad(0<r\le r_0).
\]
The uniform contraction theorem
\cite[Theorem~1.260]{Chicone2024}
therefore gives smooth dependence of $U$ on $(\beta,c)$
in $C([0,r_0],\R^N)$.

The integral formula for $rU'(r)$ gives smooth dependence of
$U'(r_0)$ as well. Ordinary differential equation theory away
from $r=0$ then gives smooth dependence of $w$ and all its
$t$ derivatives on every finite cylinder interval.

\item \emph{Uniform decay.}
At the base point, \cite[Lemma~2.1]{LZ2014} gives
\[
 \lim_{t\to\infty}w_i'(t,p_*)=-m_{*,i},
 \qquad m_{*,i}>\gamma_{*,i}.
\]
Choose $a>0$ and $T>0$ so that
\[
 \gamma_{*,i}+w_i'(T,p_*)<-3a
 \qquad(i\in I).
\]
By the finite-interval parameter dependence proved in (i),
after decreasing $\mathcal U$ we have
\[
 \gamma_i+w_i'(T,p)\le-2a
 \qquad(p\in\mathcal U,\ i\in I).
\]
Since $w_i''=-\sum_j a_{ij}q_j\le0$, this inequality persists
for every $t\ge T$. Consequently,
\[
 q_i(t,p)\le q_i(T,p)\e^{-2a(t-T)}
 \qquad(t\ge T).
\]
This is the same monotonicity argument used in
\cite[proof of Proposition~1, (23)]{LZ2014}.

At the other end, $w_i\le c_i$ gives
\[
 q_i(t,p)\le C\e^{\gamma_*t}
 \qquad(t\le0).
\]
Combining these estimates with the bounds on finite intervals,
we obtain, for some common $\nu>0$,
\[
 0<q_i(t,p)\le C\e^{-\nu|t|}
 \qquad(t\in\R,\ p\in\mathcal U).
\]
Choosing $\rho>0$ so that
$|\beta-\beta_*|+|c-c_*|<\rho$ implies $p\in\mathcal U$
proves \eqref{eq:radial-uniform-decay}.
In particular, every solution in this neighborhood has finite
component masses.

\item \emph{Parameter derivatives of the density.}
We first establish an estimate for the nonhomogeneous linear
equations arising upon differentiation. Put
\[
 B(t,p)=A\operatorname{diag}(q(t,p)).
\]
The bound in (ii) implies
\[
 \sup_{p\in\mathcal U}
 \int_0^\infty(1+s)\|B(s,p)\|\,\dd s<\infty.
\]
Differentiating the radial equation with respect to the parameters
leads to linear equations with this coefficient matrix.
Suppose that
\[
 y''+By=F,\qquad F\in L^1(0,\infty).
\]
Integrating twice and dividing by $1+t$, we find
\[
 \frac{|y(t)|}{1+t}
 \le D+\int_0^t
 (1+s)\|B(s,p)\|
 \frac{|y(s)|}{1+s}\,\dd s,
\]
where
\[
 D=|y(0)|+|y'(0)|+\int_0^\infty|F(s)|\,\dd s.
\]
Gronwall's inequality bounds the quotient on the left by $CD$,
uniformly in $p$. Integrating the equation once then gives
\[
 \sup_{t\ge0}
 \left(\frac{|y(t)|}{1+t}+|y'(t)|\right)\le CD.
\]

Let $\alpha$ be a parameter multi-index with $|\alpha|\ge1$,
and set $y_\alpha=\partial_p^\alpha w$. These derivatives
exist on finite intervals by (i). Differentiating the radial
equation gives
\[
 y_\alpha''+By_\alpha
 =-A(q_iP_{\alpha,i})_{i\in I}.
\]
For $|\alpha|=1$, we have
\[
 P_{\alpha,i}=t\,\partial_p^\alpha\beta_i.
\]
For higher orders, $P_{\alpha,i}$ is a sum of products of
parameter derivatives of $\gamma_i t+w_i$ of strictly lower
orders, with at most $|\alpha|$ factors.
The lower-order growth bounds therefore give
\[
 |q_iP_{\alpha,i}|
 \le C_\alpha(1+t)^{|\alpha|}\e^{-\nu t}
 \qquad(t\ge0).
\]
Applying the preceding linear estimate successively in
$|\alpha|$, we obtain
\[
 |y_\alpha(t,p)|\le C_\alpha(1+t),
 \qquad |y_\alpha'(t,p)|\le C_\alpha
 \qquad(t\ge0).
\]

Together with the bounded parameter derivatives of $w$ near
$t=-\infty$, these estimates control all parameter derivatives
of $q$. Repeated use of
\[
 q_i'=(\gamma_i+w_i')q_i,\qquad w''=-Aq,
\]
and the integral formula for $w'$ then gives
\begin{equation}\label{eq:density-derivatives}
 |\partial_t^j\partial_p^\alpha q(t,p)|
 \le C_{\alpha,j}(1+|t|)^{|\alpha|}
       \e^{-\nu|t|}
 \qquad(t\in\R,\ j\ge0).
\end{equation}
The neighborhood and $\nu$ remain independent of the
differentiation orders.

For every fixed $\ell\ge0$ and $0<\delta<\nu$, these estimates
bound all parameter derivatives in $C^\ell_\delta$.
Apply the pointwise Taylor formula to
$\partial_t^j\partial_p^\alpha q$, $0\le j\le\ell$,
and estimate its integral remainder by the same bounds
at parameter order $|\alpha|+2$.
Taking the weighted supremum proves smooth dependence
of $q$ in $C^\ell_\delta$.

\item \emph{Masses and asymptotic constants.}
The maps
\[
 q\longmapsto\int_{\R}q(t)\,\dd t,
 \qquad
 q\longmapsto\int_{\R}tq(t)\,\dd t
\]
are continuous linear maps on $C^0_\delta$. By (iii),
\[
 \sigma=\int_{\R}q,\qquad
 m=A\sigma,\qquad
 c^+=c+A\int_{\R}tq(t)\,\dd t
\]
therefore depend smoothly on $p$.

\item \emph{The remainder.}
The integral representation gives
\begin{equation}\label{eq:left-remainder}
 w(t)-c=-A\int_{-\infty}^t(t-s)q(s)\,\dd s,
\end{equation}
and
\begin{equation}\label{eq:right-remainder}
 w(t)+mt-c^+
 =-A\int_t^\infty(s-t)q(s)\,\dd s.
\end{equation}
Since $\chi_-+\chi_+=1$, the remainder $V$ in
\eqref{eq:weighted-smoothness} satisfies
\[
 V(t)=-A\left[
 \chi_-(t)\int_{-\infty}^t(t-s)q(s)\,\dd s
 +\chi_+(t)\int_t^\infty(s-t)q(s)\,\dd s
 \right].
\]
This is a fixed linear map of $q$. For example, when $t\le0$,
\[
 \left|\int_{-\infty}^t(t-s)q(s)\,\dd s\right|
 \le\delta^{-2}\|q\|_{C^0_\delta}\e^{\delta t}.
\]
The analogous estimate holds for the right-hand integral when
$t\ge1$. Their first derivatives satisfy the same bounds with
$\delta^{-1}$ in place of $\delta^{-2}$, and their second
derivatives equal $q$.

Derivatives of the cutoffs are supported in $[0,1]$, where
the integrals are controlled by
\[
 \int_{\R}(1+|s|)|q(s)|\,\dd s.
\]
It follows that
\[
 \|V\|_{C^{\ell+2}_\delta}
 \le C_{\ell,\delta}\|q\|_{C^\ell_\delta}.
\]
Thus the remainder also depends smoothly on $p$ in the asserted
weighted spaces.

\end{enumerate}
\end{proof}

\begin{remark}
On a fixed parameter neighborhood where $|\beta_i|<1/2$ for all
$i\in I$, the local estimates give
\[
 |\partial_p^\alpha q_i(t,p)|
 \le C_\alpha(1+|t|)^{|\alpha|}\e^{3t/2}
 \qquad\text{ for }t\le0.
\]
By \eqref{eq:radial-first},
\[
 |\partial_t\partial_p^\alpha w_i(t,p)|
 \le C_\alpha(1+|t|)^{|\alpha|}\e^{3t/2}
 \qquad\text{ for }t\le0.
\]
Consequently, since $r=\e^t$,
\[
 |\partial_r\partial_p^\alpha U_i(r,p)|
 \le C_\alpha(1+|\log r|)^{|\alpha|}r^{1/2}
 \longrightarrow0
 \qquad\text{ as }r\downarrow0.
\]
Together with the continuous extension supplied by the local
construction, this shows that the profiles and every fixed
finite-order parameter derivative extend as $C^1$ functions
of $x$ across the origin, with gradient zero there.
\end{remark}

Differentiating with respect to the initial values leads to the radial
linearized system. We determine its bounded solutions and the kernel
of the derivative of the mass map.

At fixed $\beta$, a radial solution of the linearized system regular at the origin satisfies
\begin{equation}\label{eq:radial-jacobi}
 \phi''+A\operatorname{diag}(q)\phi=0,
 \qquad \lim_{t\to-\infty}\phi(t)=h,
 \qquad \lim_{t\to-\infty}\phi'(t)=0.
\end{equation}
These conditions give a bounded continuous extension at the origin.

\begin{lemma}\label{lem:jacobi-growth}
Fix a finite-mass radial solution with parameters $(\beta,c)$ and
density $q$ as in \eqref{eq:radial-cylinder}.
For every $h\in\R^N$, the following assertions hold.
\begin{enumerate}
\renewcommand{\labelenumi}{(\roman{enumi})}

\item Problem \eqref{eq:radial-jacobi} has a unique solution given by
\[
 \phi_i(t)
 =\left.\frac{\dd}{\dd\varepsilon}
 w_i(t,\beta,c+\varepsilon h)\right|_{\varepsilon=0}
 =\sum_{j=1}^N h_j
       \frac{\partial w_i}{\partial c_j}(t,\beta,c).
\]
The integrals
\begin{equation}\label{eq:mass-variation}
 \dot\sigma_i[\phi]
 =\int_\R q_i(t)\phi_i(t)\,\dd t
\end{equation}
are absolutely convergent.

\item For some $\delta>0$ and $k\in\R^N$,
\begin{align}
 \phi(t)&=h+O(\e^{\delta t})
 &&\text{as }t\to-\infty,
 \label{eq:jacobi-left}\\
 \phi(t)&=-A\dot\sigma[\phi]\,t+k+O(\e^{-\delta t})
 &&\text{as }t\to+\infty.
 \label{eq:jacobi-asymptotic}
\end{align}
Consequently, $\phi$ is bounded on $\R$ if and only if
$\dot\sigma[\phi]=0$.

\end{enumerate}
\end{lemma}

\begin{proof}
\begin{enumerate}
\renewcommand{\labelenumi}{(\roman{enumi})}

\item \emph{Initial-value and mass variations.}
Proposition~\ref{prop:radial-parameters} shows that the derivative
$\phi$ is well-defined. Differentiating \eqref{eq:radial-ode}
at fixed $\beta$ gives the differential equation in
\eqref{eq:radial-jacobi}. The left-end parameter estimates in
that proposition also give
\[
 \phi(t)\longrightarrow h,\qquad
 \phi'(t)\longrightarrow0
 \qquad(t\to-\infty),
\]
so $\phi$ satisfies the prescribed limits.

For uniqueness, the difference of two solutions of
\eqref{eq:radial-jacobi}, viewed as a radial function of
$r=\e^t$, extends continuously to the origin with value zero.
It therefore vanishes by \cite[Lemma~2.5]{LZ2014}.

The parameter estimates in
Proposition~\ref{prop:radial-parameters} justify differentiation
under the mass integral and give
\[
\begin{aligned}
 \left.\frac{\dd}{\dd\varepsilon}
 \sigma_i(\beta,c+\varepsilon h)\right|_{\varepsilon=0}
 =
 \int_\R
 \left.\frac{\dd}{\dd\varepsilon}
 q_i(t,\beta,c+\varepsilon h)\right|_{\varepsilon=0}\,\dd t=\int_\R q_i(t)\phi_i(t)\,\dd t
 =\dot\sigma_i[\phi].
\end{aligned}
\]
In particular, the integral in \eqref{eq:mass-variation}
is absolutely convergent.

\item \emph{Asymptotics and boundedness.}
Set
\[
 k=\left.\frac{\dd}{\dd\varepsilon}
 c^+(\beta,c+\varepsilon h)\right|_{\varepsilon=0}.
\]
For sufficiently small $\delta>0$, differentiate the remainder
\eqref{eq:weighted-smoothness} along the same parameter curve.
Its smooth dependence in $C^1_\delta$, together with
$m=A\sigma$ and the mass-variation identity proved in (i), gives
\[
 \phi-h\chi_-
 -(k-A\dot\sigma[\phi]\,t)\chi_+
 \in C^1_\delta.
\]
For $t\le0$, this says that $\phi-h$ and $\phi'$ are
$O(\e^{\delta t})$, proving \eqref{eq:jacobi-left}.
For $t\ge1$, it gives
\[
 \phi(t)=-A\dot\sigma[\phi]\,t+k+O(\e^{-\delta t}),
\]
which is \eqref{eq:jacobi-asymptotic}.

Finally, $\phi$ is bounded at the left end, and
\eqref{eq:jacobi-asymptotic} shows that it is bounded at the
right end if and only if $A\dot\sigma[\phi]=0$.
Since $A$ is invertible, this is equivalent to $\dot\sigma[\phi]=0$.

\end{enumerate}
\end{proof}

\subsection{Dependence on the component masses}

Theorem \ref{thm:published-regular-kernel} (ii) allows us to replace the initial values, modulo dilation, by the component masses. The mass surface can also be described locally.

\begin{lemma}\label{lem:mass-graph}
For every $\sigma\in\Gamma_\beta$, we have $d_i=(A\sigma)_i-\gamma_i>0$ for all $i\in I$. Near every pair $(\beta_*,\sigma_*)$ with $\sigma_*\in\Gamma_{\beta_*}$, the condition $\sigma\in\Gamma_\beta$ is equivalent to $\sigma_N=g(\beta,\sigma_1,\ldots,\sigma_{N-1})$ for a smooth function $g$.
\end{lemma}
\begin{remark} 
In follows, we set $z=(\sigma_1,\ldots,\sigma_{N-1})$, we write
\begin{equation}\label{eq:mass-graph}
 \Sigma(\beta,z):=(z_1,\ldots,z_{N-1},g(\beta,z)).
\end{equation}
\end{remark}
\begin{proof}[Proof of Lemma \ref{lem:mass-graph}]
For every $i\in I$, the defining conditions of $\Gamma_\beta$ give
\begin{align}
 0<\Lambda_{I\setminus\{i\}}^\beta(\sigma)
 = -2\gamma_i\sigma_i
    +2\sigma_i\sum_{j\ne i}a_{ij}\sigma_j
    +a_{ii}\sigma_i^2=2\sigma_i d_i-a_{ii}\sigma_i^2.\nonumber
\end{align}
Since $\sigma_i>0$ and $a_{ii}\ge0$, we obtain $d_i>0$.
By symmetry of $A$,
\[
 \frac{\partial\Lambda_I^\beta}{\partial\sigma_i}(\sigma)
 =2\gamma_i-2(A\sigma)_i=-2d_i.
\]
In particular, the derivative with respect to $\sigma_N$ is
nonzero at $(\beta_*,\sigma_*)$. The implicit function theorem
therefore gives the existence of the smooth function $g$.
\end{proof}

At a finite-mass radial profile with exponents $\beta$ and mass $\sigma$, let
\[
 \mathcal E_{\beta,\sigma}
 =\{\phi:\ \phi\text{ is a radial solution of }\eqref{eq:radial-jacobi} \text{ for some }h\in\R^N\}.
\]
Thus its elements have finite limits $h$ and derivative tending to zero as $t\to-\infty$. No boundedness as $t\to+\infty$ is imposed.

In follows, we write $d^\perp:=\{\xi\in\R^N:d\cdot\xi=0\}$.

\begin{proposition}\label{thm:mass-chart}
Fix a finite-mass radial solution at $(\beta_*,\sigma_*)$.
For any couple $(\beta,z)$ sufficiently close to $(\beta_*,z_*)$ where $z_*=(\sigma_{*,1},...\sigma_{*,N-1})$ there exists a unique radial solution
$\widetilde U_{\beta,z}$ satisfying
\[
 \sigma(\widetilde U_{\beta,z})=\Sigma(\beta,z),
 \qquad \widetilde U_N(0)=0.
\]
For every fixed $\ell\ge0$ and $0<\delta<\nu$ as in Proposition~\ref{prop:radial-parameters}, the map
\begin{align}\label{def:mapbetaz}
 (\beta,z)\longmapsto
 \bigl(c,c^+,q,w-c\chi_--(c^+-mt)\chi_+\bigr)
 \in\R^{2N}\times(C^\ell_\delta)^2
\end{align}
is $C^\infty$. Moreover, with $\dot\sigma:\mathcal E_{\beta,\sigma}\to \R^N$ defined by \eqref{eq:mass-variation}, we have
\[
 \Ran\dot\sigma=d^\perp=(A\sigma-\gamma)^\perp
 \qquad\text{on }\mathcal E_{\beta,\sigma},
 \qquad \sigma=\Sigma(\beta,z).
\]
\end{proposition}

\begin{proof}
Fix $\beta$. By Lemma~\ref{lem:jacobi-growth}, for every $h\in\R^N$ there exists a unique solution $\phi\in\mathcal E_{\beta,\sigma}$ satisfying $\lim_{t\to-\infty}\phi(t)=h$. Differentiating $\Lambda_I^\beta(\sigma)=0$ along the initial-value family gives
\[
\begin{aligned}
 0=\left.\frac{\mathrm d}{\mathrm d\varepsilon}\Lambda_I^\beta\bigl(\sigma(\beta,c+\varepsilon h)\bigr)\right|_{\varepsilon=0}=\sum_i\left(2\gamma_i-2\sum_j a_{ij}\sigma_j\right)\dot\sigma_i[\phi]=-2\sum_i d_i\dot\sigma_i[\phi].
\end{aligned}
\]
Thus the image is contained in $d^\perp$. The same lemma shows that $\dot\sigma[\phi]=0$ exactly when $\phi$ is bounded. This is due to Theorem~\ref{thm:published-regular-kernel} (ii), this is equivalent to $\phi\in\Span\{Z\}$. Since $\mathcal E_{\beta,\sigma}$ has dimension $N$,
\begin{equation}\label{eq:mass-variation-kernel-range}
 \Ran\dot\sigma=d^\perp=(A\sigma-\gamma)^\perp,\qquad \Ker\dot\sigma=\Span\{Z\}.
\end{equation}

The normalization $c_N=0$ fixes dilation uniquely because
$\gamma_N>0$. Write $c'=(c_1,\ldots,c_{N-1})$ and
$\sigma'=(\sigma_1,\ldots,\sigma_{N-1})$. If a vector $h=(h',0)$
satisfies $D_{c'}\sigma'h'=0$, the preceding identity and $d_N>0$
give $\dot\sigma[\phi]=0$. Hence $\phi=aZ$, and its initial limit
satisfies $0=h_N=a\gamma_N$, so $h'=0$.
Consequently $D_{c'}\sigma'$ is invertible. The smooth map
\[
 (\beta,c')\longmapsto(\beta,\sigma'(\beta,c',0))
\]
has an invertible derivative and the inverse function theorem gives
smooth initial values $c'(\beta,z)$. The last mass is $g(\beta,z)$ by Lemma~\ref{lem:mass-graph}. Composition with Proposition~\ref{prop:radial-parameters} gives the smoothness of (\ref{def:mapbetaz}). Theorem~\ref{thm:published-classification} and the
normalization give uniqueness.
\end{proof}

\begin{remark}\label{rem:fixed-mass}
At fixed $\beta$, differentiating
$\Lambda_I^\beta(\Sigma(\beta,z))=0$ gives
\[
 0=-2d_l-2d_N\partial_{z_l}g.
\]
Hence
\begin{equation}\label{eq:g-derivatives}
 \partial_{z_l}g=-\frac{d_l}{d_N},
 \qquad 1\le l<N,\qquad \sigma=\Sigma(\beta,z).
\end{equation}
Differentiating \eqref{eq:mass-graph} therefore yields
\[
 \partial_{z_l}\Sigma
 =e_l+(\partial_{z_l}g)e_N
 =e_l-\frac{d_l}{d_N}e_N
 \qquad(1\le l<N).
\]

If the same vector $\sigma$ belongs to both $\Gamma_\beta$
and $\Gamma_{\beta_*}$, then
\[
 0=\Lambda_I^\beta(\sigma)-\Lambda_I^{\beta_*}(\sigma)
  =2\sigma\cdot(\beta-\beta_*).
\]
Hence
\begin{equation}\label{eq:fixed-mass-hyperplane}
 \sigma\cdot(\beta-\beta_*)=0.
\end{equation}

For every nonzero vector $\xi\in d^\perp$
at fixed $\beta$, Proposition~\ref{thm:mass-chart} gives a
corresponding $\phi\in\mathcal E_{\beta,\sigma}$.
Equation~\eqref{eq:jacobi-asymptotic} gives
\[
 \phi(t)=-A\xi\,t+O(1)
 \qquad(t\to+\infty).
\]
Since $A$ is invertible, $A\xi\ne0$, so $\phi$ is unbounded
in at least one component.
\end{remark}

\subsection{Dilation}

The masses do not determine the dilation parameter. The median condition
\eqref{eq:median} selects a unique representative and varies smoothly
with the radial family. It replaces the initial-value normalization
$\widetilde U_N(0)=0$.

\begin{proposition}\label{prop:median}
Let $w$ correspond to a radial solution with finite component masses
and parameters $p=(\beta,c)$. For $s\in\R$, define
\begin{equation}\label{eq:dilation-action}
 (T_sw)_i(t)=w_i(t+s)+\gamma_i s,
 \qquad i\in I.
\end{equation}
There exists a unique $s=s(p)$ such that $T_sw$ satisfies
\eqref{eq:median}. Equivalently,
\begin{equation}\label{eq:median-shift}
 \sum_i\int_{-\infty}^{s(p)}q_i(t,p)\dd t
 =\frac12\sum_i\sigma_i(p).
\end{equation}
Moreover, the function $s(p)$ is locally smooth and the corresponding solutions have the smooth dependence stated in Proposition~\ref{prop:radial-parameters}. 
\end{proposition}
\begin{proof}
Set
\[
 H(p,s):=\sum_i\int_{-\infty}^sq_i(t,p)\dd t
       -\frac12\sum_i\sigma_i(p).
\]
Proposition~\ref{prop:radial-parameters} makes $H$ smooth.
For fixed $p$, its limits at $s=\pm\infty$ are respectively
$\pm\frac12\sum_i\sigma_i(p)$ and
$\partial_s H=\sum_iq_i(s,p)>0$. Thus it has a unique zero, smooth in $p$
by the implicit function theorem.

The solution $T_{s(p)}w$ has initial values $c+\gamma s(p)$. By uniqueness in Proposition~\ref{prop:radial-parameters}(i),
\[
 (T_{s(p)}w)(t,p)
 =w\bigl(t,\beta,c+\gamma s(p)\bigr).
\]
The map $p\mapsto(\beta,c+\gamma s(p))$ is smooth. Moreover, dilation preserves the component masses, so the solutions on the right have finite component masses. Applying Proposition~\ref{prop:radial-parameters} near each such parameter pair gives all the asserted smoothness.
\end{proof}

\subsection{Uniform estimates and condition (\ref{eq:NC})}
To begin with, let us introduce the following definition.
\begin{definition}
Let $\widetilde U_{\beta,z}$ be the solution in
Proposition~\ref{thm:mass-chart}. We denote by $U_{\beta,z}$
its unique dilation satisfying \eqref{eq:median}.
\end{definition}

We next establish uniform estimates for these solutions when
$\beta=0$ and their mass vectors range over a compact subset of
$\Gamma_0$. These estimates will be used to construct vectors $b$
with both positive and negative components satisfying \eqref{eq:NC}.

For $\eta\in\Gamma_0$, recall the $U^\eta$ is the median-normalized regular radial solution of mass $\eta$ and we define
\[
 h_i^\eta=\e^{U_i^\eta(0)},
 \qquad
 E_i^\eta(r)=\int_{B_r}\e^{U_i^\eta(x)}\dd x.
\]
The latter is introduced in (\ref{eq:disk-integral}).

\begin{lemma}\label{lem:compact-profiles}
Let $K\subset\Gamma_0$ be nonempty and compact. All estimates below and in the proof are uniform for $\eta\in K$. Their constants may depend on $K$ and the fixed orders of differentiation.
\begin{enumerate}
\renewcommand{\labelenumi}{(\roman{enumi})}

\item The initial values $c$ and constants $c^+$ are uniformly bounded.
On a finite cover of $K$, the parameter and $t$ derivative estimates
above are uniform at every fixed order, with a common positive
exponential decay rate.

\item As $r\downarrow0$,
\begin{align}
 U_i^\eta(r)
 &=U_i^\eta(0)-\frac{r^2}{4}\sum_j a_{ij}h_j^\eta+O(r^4),
 \label{eq:regular-expansion}\\
 E_i^\eta(r)
 &=\pi h_i^\eta r^2+O(r^4).
 \label{eq:disk-integral-zero}
\end{align}
These expansions may be differentiated to any fixed order in $r$
and in the local parameters on $\Gamma_0$, with the corresponding
powers of $r$ in the remainders.

\item There exists $\epsilon>0$, depending on $K$, such that
\begin{equation}\label{eq:disk-integral-infty}
 E_i^\eta(r)=2\pi\eta_i+O(r^{-\epsilon})
 \qquad(r\to\infty).
\end{equation}

\end{enumerate}
\end{lemma}

\begin{proof}
\begin{enumerate}
\renewcommand{\labelenumi}{(\roman{enumi})}

\item Proposition~\ref{thm:mass-chart} and
Proposition~\ref{prop:median} give solutions depending smoothly
on the parameters near every $\eta\in\Gamma_0$.
On overlapping neighborhoods, these solutions agree by
Theorem~\ref{thm:published-classification} and
condition \eqref{eq:median}.
By compactness of $K$, the local estimates above are uniform on
a finite cover of $K$, with a common positive decay rate and
constants depending on the fixed differentiation orders.
In particular,
\begin{equation}\label{eq:h-bounds}
 0<h:=\min_{\substack{\eta\in K\\i\in I}}h_i^\eta
 \le\max_{\substack{\eta\in K\\i\in I}}h_i^\eta<\infty.
\end{equation}

\item At $\beta=0$, the radial integral equation gives
\[
 U_i^\eta(r)-U_i^\eta(0)
 =-\sum_j a_{ij}\int_0^r
       \log\frac rs\,s\e^{U_j^\eta(s)}\,\dd s.
\]
The bounds for the initial values and monotonicity first give
$U_i^\eta(r)-U_i^\eta(0)=O(r^2)$, uniformly for $\eta\in K$.
Hence for small $s$
\[
 \e^{U_j^\eta(s)}=h_j^\eta+O(s^2).
\]
Using
\[
 \int_0^r\log\frac rs\,s\,\dd s=\frac{r^2}{4},
 \qquad
 \int_0^r\log\frac rs\,s^3\,\dd s=\frac{r^4}{16},
\]
we obtain \eqref{eq:regular-expansion}.
Integrating $2\pi r\e^{U_i^\eta(r)}$ then gives
\eqref{eq:disk-integral-zero}.

Repeated use of
\[
 (U_i^\eta)'(r)
 =-\frac1r\sum_j a_{ij}\int_0^r
                      s\e^{U_j^\eta(s)}\,\dd s
\]
and the radial equation gives the differentiated remainder
estimates in $r$.

At fixed $\beta=0$, Propositions~\ref{prop:radial-parameters}
and~\ref{thm:mass-chart} give
\begin{equation}\label{eq:regular-sharp-left}
 \partial_\eta^{\mathfrak a}\partial_t^j w'(t)
 =O(\e^{2t})
 \qquad(t\to-\infty)
\end{equation}
for every fixed parameter order and every fixed $j\ge0$.
These estimates are preserved by the normalization in
Proposition~\ref{prop:median} and give the asserted uniform
remainder estimates through the radial integral formulas.

\item By the definition of the component masses,
\[
 2\pi\eta_i-E_i^\eta(r)
 =2\pi\int_{\log r}^\infty q_i^\eta(t)\,\dd t.
\]
The common exponential bound from (i) and
\eqref{eq:radial-uniform-decay} imply
\eqref{eq:disk-integral-infty}, with a fixed positive
$\epsilon$ depending on $K$.

\end{enumerate}
\end{proof}

\begin{proposition}\label{prop:nc-directions}
For every nonempty compact $K\subset\Gamma_0$, there exist a nonempty compact set
$B\subset\{b\in\R^N:\sum_i|b_i|=1\}$ and a number $\kappa>0$ such that
every $b\in B$ has at least one positive and one negative component and
\begin{equation}\label{eq:nc-uniform}
 \left|\sum_i b_iE_i^\eta(r)\right|
 \ge\kappa\sum_i|b_i|E_i^\eta(r)
 \qquad(\eta\in K,\ b\in B,\ r>0).
\end{equation}
\end{proposition}

\begin{proof}
We construct $b$ with $b_1>0$, $b_2<0$ and $b_i=0$ for $i\ge3$.
The expansions \eqref{eq:disk-integral-zero} and
\eqref{eq:disk-integral-infty}, together with positivity on compact
radius intervals, give
\[
 C=\sup_{\eta\in K,\ r>0}\frac{E_2^\eta(r)}{E_1^\eta(r)}<\infty.
\]
Choose $0<a<C^{-1}$ and set
$b_0=(1,-a,0,\ldots,0)/(1+a)$. Then
\[
 \frac{\sum_i (b_0)_iE_i^\eta(r)}{\sum_i|(b_0)_i|E_i^\eta(r)}
 =\frac{E_1^\eta(r)-aE_2^\eta(r)}{E_1^\eta(r)+aE_2^\eta(r)}
 \ge\frac{1-aC}{1+aC}>0.
\]
Taking $B=\{b_0\}$ and $\kappa=(1-aC)/(1+aC)$ proves the assertion.
\end{proof}

\section{Integral identities}
\label{sec:moments}

Assume that $A$ satisfies \textup{(H1)}.
We first prove an identity for a possibly nonradial singular solution.
We then compute the inverse moment of a translated regular radial
solution, for an arbitrary signed vector of coefficients. Through this section, we denote
\[
 f_i(x)=|x|^{\beta_i}\e^{u_i(x)},\qquad
 \mathcal K(x)=\frac{x}{|x|^2}\quad(x\ne0).
\]
The kernel $\mathcal K$ is locally integrable in two dimensions.  The
integrals below involve an inverse power of the distance to the singular
point. Their convergence does not follow from finite mass alone.

\begin{lemma}\label{mom:origin}
Let $u$ be a finite-mass solution of \eqref{eq:system}, with each $u_i$
bounded near the origin and assume $\beta_i>-1$ for every $i$.
Then each $u_i$ extends as a $C^1$ function across the origin.  For every $i\in I$,
\[
 \int_{\R^2}\frac{f_i(x)}{|x|}\,\dd x<\infty,
 \qquad
 \int_{\R^2}f_i(x)|\nabla u_i(x)|\,\dd x<\infty.
\]
Moreover, $|\nabla u_i|$ is bounded near zero and is $O(|x|^{-1})$
as $|x|\to\infty$.
\end{lemma}

\begin{proof}
Choose $R>0$ such that every $u_i$ is bounded in
$B_R\setminus\{0\}$. Choose $p>2$ sufficiently close to $2$ that $p\beta_i>-2$ for every
negative $\beta_i$.  Boundedness of $u_i$ near zero gives $f_i\in L^p(B_R)$ for a sufficiently small disk.  Consequently,
\[
 g_i=\sum_j a_{ij}f_j\in L^p(B_R).
\]
Define the local Newton potential
\[
 v_i(x)=-\frac1{2\pi}\int_{B_R}\log|x-y|\,g_i(y)\,\dd y.
\]
On a smaller disk this potential and its gradient are bounded.  Indeed,
if $p'=p/(p-1)$, then $p'<2$, and both the logarithmic kernel and the
kernel $|x-y|^{-1}$ have locally bounded $L^{p'}$ norms, uniformly in
$x$ in that smaller disk.  H\"older's inequality proves the bounds.
Approximating $g_i$ in $L^p$ by smooth functions proves continuity of
the gradient as well.  Since $-\Delta v_i=g_i$, the difference
$u_i-v_i$ is bounded and harmonic on a punctured disk.  The removable
singularity theorem for bounded harmonic functions extends it
harmonically across zero.  This proves the asserted $C^1$ extension
and the local gradient bound.

The asymptotic estimates in Theorem~\ref{thm:published-asymptotics}
give, with $m=A\sigma$ and $d_i=m_i-2-\beta_i>0$,
\[
 f_i(x)\le C|x|^{-2-d_i},\qquad
 |\nabla u_i(x)|\le C|x|^{-1}
 \quad (|x|\ge R_0).
\]
These constants are allowed to depend on the solution.  Near the origin,
\[
 \int_{B_r}\frac{f_i(x)}{|x|}\,\dd x
 \le C\int_0^r s^{\beta_i}\,\dd s
 =\frac{C}{1+\beta_i}r^{1+\beta_i}.
\]
On $\R^2\setminus B_R$, we have
\[
 \int_{\R^2\setminus B_R}\frac{f_i(x)}{|x|}\,\dd x
 \le \frac1R\int_{\R^2\setminus B_R}f_i(x)\,\dd x
 <\infty.
\]
Thus the inverse moment is absolutely convergent. The local gradient bound, the decay estimate for $|x|\ge R_0$ and integrability of $f_i$ imply
\[
 \int_{\R^2}f_i(x)|\nabla u_i(x)|\,\dd x<\infty.
\]
\end{proof}

\begin{lemma}\label{cmp:gradient}
Under the assumptions of Lemma~\ref{mom:origin},
\begin{equation}\label{cmp:gradient-representation}
 \nabla u_i(x)=-\frac1{2\pi}\sum_j a_{ij}
 \int_{\R^2}\mathcal K(x-y)f_j(y)\,\dd y.
\end{equation}
The integral is absolutely convergent at every $x\in\R^2$.
\end{lemma}

\begin{proof}
Choose $p>2$ as in Lemma~\ref{mom:origin}, whose proof gives
$f_j\in L^p(B_R)$ near the origin. Continuity away from the origin and the decay in Theorem~\ref{thm:published-asymptotics} then give $f_j\in L^p(\R^2)$.
The finite-mass assumption gives $f_j\in L^1(\R^2)$.
Since $p'=p/(p-1)<2$,
\[
 \sup_x\int_{\R^2}|\mathcal K(x-y)|f_j(y)\,\dd y
 \le C_p\|f_j\|_{L^p}+\|f_j\|_{L^1}<\infty.
\]
On $|x-y|<1$ we use H\"older's inequality. On $|x-y|\geq 1$ we use $|\mathcal K|\le1$. These estimates also justify differentiation of
\eqref{eq:published-potential}. Indeed, on a fixed compact set, approximation
in local $L^p$ handles the kernel near $x$, and dominated convergence
applies away from $x$ and at infinity. The terms independent of $x$
differentiate to zero, giving \eqref{cmp:gradient-representation}.
\end{proof}

The next identity retains the signs of the exponents. Symmetry of $A$
and oddness of $\mathcal K$ cancel the interaction terms after summation.

\begin{proposition}\label{mom:identity}
Under the assumptions of Lemma~\ref{mom:origin},
\begin{equation}\label{mom:identity-eq}
 \sum_{i=1}^N\beta_i
 \int_{\R^2}\mathcal K(x)f_i(x)\,\dd x=0.
\end{equation}
\end{proposition}
\begin{proof}
On $\R^2\setminus\{0\}$, differentiating $f_i=|x|^{\beta_i}\e^{u_i}$
gives
\[
 \nabla f_i=\beta_i\mathcal K f_i+f_i\nabla u_i.
\]
Both terms on the right are integrable by Lemma~\ref{mom:origin}.
The boundary terms on $\partial B_r$ and $\partial B_R$ satisfy
\[
 \left|\int_{\partial B_r}f_i\frac{x}{r}\,\dd S\right|
 \le Cr^{1+\beta_i},
 \qquad
 \left|\int_{\partial B_R}f_i\frac{x}{R}\,\dd S\right|
 \le CR^{-1-d_i}.
\]
Since $\beta_i>-1$ and $d_i>0$, both terms tend to zero as
$r\downarrow0$ and $R\to\infty$.
Integrating the differential identity over $B_R\setminus\overline{B_r}$
and passing to these limits therefore gives
\[
 \beta_i\int_{\R^2}\mathcal K(x)f_i(x)\,\dd x
 +\int_{\R^2}f_i(x)\nabla u_i(x)\,\dd x=0.
\]

We now compute the second term after summation over $i$.
The estimate in the proof of Lemma~\ref{cmp:gradient} gives
\[
 \iint_{\R^2\times\R^2}
 |\mathcal K(x-y)|f_i(x)f_j(y)\,\dd x\,\dd y
 \le \|f_i\|_{L^1(\R^2)}
 \bigl(C_p\|f_j\|_{L^p(\R^2)}+\|f_j\|_{L^1(\R^2)}\bigr)
 <\infty.
\]
Thus we may substitute \eqref{cmp:gradient-representation}
and interchange the integrals.
Interchanging $(i,x)$ with $(j,y)$ and using $a_{ij}=a_{ji}$ yields
\begin{equation}\label{mom:summed-interaction-zero}
 \begin{aligned}
 \sum_i\int_{\R^2}f_i\nabla u_i
 &=-\frac1{2\pi}\sum_{i,j}a_{ij}
 \iint_{\R^2\times\R^2}
 \mathcal K(x-y)f_i(x)f_j(y)\,\dd x\,\dd y\\
 &=-\frac1{4\pi}\sum_{i,j}a_{ij}
 \iint_{\R^2\times\R^2}
 [\mathcal K(x-y)+\mathcal K(y-x)]
 f_i(x)f_j(y)\,\dd x\,\dd y\\
 &=0.
 \end{aligned}
\end{equation}
The last equality follows from the oddness of $\mathcal K$.
Summing the integrated differential identity over $i$ and applying
\eqref{mom:summed-interaction-zero} proves the proposition.
\end{proof}

\begin{remark}
The proposition gives one summed vector identity.  In general, it does
not give an identity with zero right-hand side for each component
separately.  Its range $\beta_i>-1$ is also more restrictive than the
range $\beta_i>-2$ used for the radial mass classification.
\end{remark}

\begin{proposition}\label{mom:center}
Let $U$ be a radial finite-mass solution of \eqref{eq:system}
with $\beta=0$ and centered at the origin. Set
\[
 \rho_i(x)=\e^{U_i(x)},\qquad
 E_i(r)=\int_{B_r}\rho_i(x)\,\dd x,\qquad
 \eta_i=\frac1{2\pi}\int_{\R^2}\rho_i.
\]
For $b=(b_1,\dots,b_N)\in\R^N$, we define
\begin{align}\label{def Sb}
 S_b(r)=\sum_i b_i E_i(r),\qquad
 M_b(a)=\sum_i b_i
 \int_{\R^2}\mathcal K(x)\rho_i(x-a)\,\dd x.
\end{align}
These integrals are absolutely convergent and
\begin{equation}\label{mom:center-identity}
 M_b(0)=0,\qquad
 M_b(a)=\frac{a}{|a|^2}S_b(|a|)\quad(a\ne0).
\end{equation}
\end{proposition}

\begin{proof}
Integrating the radial equation of $U$, we obtain
\[
 rU_i'(r)=-\sum_j a_{ij}\int_0^r s\e^{U_j(s)}\,\dd s\le0,
\]
so $\rho_i$ is bounded. Splitting the integral into $B_1$ and
its complement gives
\[
 \int_{\R^2}|\mathcal K(x)|\rho_i(x-a)\,\dd x
 \le 2\pi\|\rho_i\|_{L^\infty(\R^2)}
      +\|\rho_i\|_{L^1(\R^2)}
 <\infty
 \qquad(a\in\R^2).
\]
Thus the integrals defining $M_b(a)$ are absolutely convergent. Writing $F_b=\sum_i b_i\rho_i$ and changing variables yield
\[
 M_b(a)=\int_{\R^2}\mathcal K(a+y)F_b(y)\,\dd y.
\]
For $r=|a|>0$ and $s\ne r$, we get
\[
\int_0^{2\pi}\log|a+s(\cos\theta,\sin\theta)|\,\dd\theta
 =2\pi\log\max\{r,s\}.
\]
Differentiation with respect to $a$, away from
$r=s$, gives
\[
 \int_0^{2\pi}\mathcal K(a+s(\cos\theta,\sin\theta))\,\dd\theta
 =\begin{cases}2\pi a/r^2,&s<r,\\0,&s>r.\end{cases}
\]
We may therefore integrate in polar coordinates. Hence
\[
 M_b(a)=\frac{2\pi a}{r^2}\int_0^r sF_b(s)\,\dd s
       =\frac{a}{r^2}S_b(r).
\]
For $a=0$, oddness of $\mathcal K$ and radial symmetry of $F_b$ give
$M_b(0)=0$.
\end{proof}

\begin{corollary}\label{mom:nc-consequences}
In the setting of Proposition~\ref{mom:center}, let $b\in\R^N$
satisfy $\sum_i|b_i|=1$.
Assume that there exists $\kappa>0$ such that
\begin{equation}\label{mom:fixed-nc}
 |S_b(r)|\ge\kappa\sum_i|b_i|E_i(r)
 \qquad\text{for every }r>0.
\end{equation}
Here, $S_b(r)$ is defined as in (\ref{def Sb}). Then $M_b(a)$ has no zero with $a\ne0$ and
\begin{equation}\label{mom:fixed-nc-gaps}
 \left|\sum_i b_i\e^{U_i(0)}\right|
 \ge\kappa\sum_i|b_i|\e^{U_i(0)},\qquad
 |b\cdot\eta|\ge\kappa\sum_i|b_i|\eta_i>0.
\end{equation}
For the normalized family $U^\eta$, $\eta\in K$, and a compact
set of vectors $B$ satisfying \eqref{eq:nc-uniform}, set
$c_b^\eta=\sum_i b_ih_i^\eta$ and
$h=\min_{\eta\in K,i\in I}h_i^\eta$. Then $h>0$ and
\begin{equation}\label{eq:central-gap}
 |c_b^\eta|\ge\kappa h>0
 \qquad((\eta,b)\in K\times B).
\end{equation}
\end{corollary}

\begin{proof}
Notice that $E_i(r)>0$ when $r>0$, so the right side of \eqref{mom:fixed-nc} is positive. Identity
\eqref{mom:center-identity} gives (\ref{mom:fixed-nc}). Dividing
\eqref{mom:fixed-nc} by $\pi r^2$ and letting $r\downarrow0$ gives
the first inequality. Letting $r\to\infty$ gives (\ref{mom:fixed-nc-gaps}). For the normalized family, Lemma~\ref{lem:compact-profiles} gives
$h_i^\eta=\e^{U_i^\eta(0)}>0$ continuously on $K$.
Compactness and $\sum_i|b_i|=1$ prove \eqref{eq:central-gap}.
\end{proof}

\begin{remark}\label{rem:nc-scope}
A common dilation replaces $E_i^\eta(r)$ by $E_i^\eta(\lambda r)$.
Hence condition \eqref{eq:nc-uniform} and the quantity
\[
 \sup_{\eta\in K,\ r>0}\frac{E_2^\eta(r)}{E_1^\eta(r)}
\]
remain the same under dilation. The uniform positive lower bound for $h_i^\eta$ uses the normalization \eqref{eq:median}.

Condition \eqref{eq:nc-uniform} concerns solutions with $\beta=0$.
For a radial solution $U$ with $\beta_i>-2$, continuity at the
origin gives
\[
 \int_{B_r}|x|^{\beta_i}\e^{U_i(x)}\,\dd x
 =\frac{2\pi\e^{U_i(0)}}{2+\beta_i}
  r^{2+\beta_i}(1+o(1))
 \qquad(r\downarrow0).
\]
Consequently,
\[
 \frac{\displaystyle\int_{B_r}|x|^{\beta_i}\e^{U_i(x)}\,\dd x}
      {\displaystyle\int_{B_r}|x|^{\beta_j}\e^{U_j(x)}\,\dd x}
 =
 \frac{(2+\beta_j)\e^{U_i(0)}}{(2+\beta_i)\e^{U_j(0)}}
 r^{\beta_i-\beta_j}(1+o(1)).
\]

The function $c_b^\eta=\sum_i b_ih_i^\eta$ is continuous in
$(\eta,b)$.
By compactness of $K\times B$ and \eqref{eq:central-gap},
\[
 |c_b^\eta|\ge\frac{\kappa h}{2}
\]
for all $b\in B$ and all $\eta\in\Gamma_0$ sufficiently close to $K$.
Condition \eqref{eq:nc-uniform} is still assumed only for
$\eta\in K$.
\end{remark}

\begin{example}\label{mom:mixed-example}
The hypotheses permit an indefinite coupling matrix.  Take
\[
 A=\begin{pmatrix}0&1\\1&0\end{pmatrix},
 \qquad
 U_1(r)=U_2(r)=\log 8-2\log(1+r^2).
\]
The matrix is symmetric, entrywise nonnegative, irreducible, and invertible,
with eigenvalues $1$ and $-1$.  Moreover,
\[
 -\Delta U_i=\frac{8}{(1+r^2)^2}
 =\sum_{j=1}^2 a_{ij}\e^{U_j},
\]
and direct radial integration gives
\[
 E_1(r)=E_2(r)
 =16\pi\int_0^r\frac{s}{(1+s^2)^2}\,\dd s
 =8\pi\frac{r^2}{1+r^2},
 \qquad \eta_1=\eta_2=4.
\]
Thus $\eta=(4,4)\in\Gamma_0$, since
\[
 \Lambda_I^0(\eta)=0,\qquad
 \Lambda_{\{i\}}^0(\eta)=16>0
 \quad(i=1,2).
\]
The vector $b=(2/3,-1/3)$ satisfies
\eqref{mom:fixed-nc} with $\kappa=1/3$ at every radius.
In contrast, for $b=(1/2,-1/2)$ the integral $M_b$ vanishes
identically.  Nonzero $b$ alone therefore gives no invertibility
conclusion for this integral.
\end{example}

\section{Compactness for small signed exponents}
\label{sec:compactness}

We prove compactness for small signed exponents up to a dilation about the
singular point. The argument first excludes point masses of the limiting measures and then recovers every component and its full mass. A point mass at $p$ means that $\mu(\{p\})>0$ for the finite nonnegative measure $\mu$ under consideration.

For a solution with positive total mass, the mapping
\[
 R\longmapsto\sum_i\int_{B_R}|x|^{\beta_i}\e^{u_i(x)}\,\dd x
\]
is continuous and strictly increasing. Thus there exists a unique $R>0$ such that
\[
 \sum_i\int_{B_R}f_i(x)\,\dd x
 =\frac12\sum_i\int_{\R^2}f_i(x)\,\dd x.
\]
The dilation
\begin{equation}\label{cmp:median-dilation}
 \widetilde u_i(y)=u_i(Ry)+(2+\beta_i)\log R,
 \qquad
 \widetilde f_i(y)=R^2f_i(Ry)
\end{equation}
preserves the equation and the masses componentwisely.  The dilated solution satisfies \eqref{eq:median} and the singular point remains at the origin.
The task of this section is to show that

\begin{theorem}\label{cmp:compactness}
Let $u_n$ be solutions of \eqref{eq:system}--\eqref{eq:mass} satisfying
\[
 \begin{gathered}
 \beta_n\to0,\qquad \sigma_n\to\sigma\in\Gamma_0,\qquad
 \beta_n=\tau_n b_n,\qquad \tau_n>0,\qquad \tau_n\to0,\\
 \sum_i|b_{i,n}|=1,\qquad b_n\to b.
 \end{gathered}
\]
\begin{equation}\label{cmp:median}
 \sum_i\int_{B_1} f_{i,n}\,\dd x
 =\pi\sum_i\sigma_{i,n},\qquad
 f_{i,n}=|x|^{\beta_{i,n}}\e^{u_{i,n}}.
\end{equation}
If $b\cdot\sigma\ne0$, a subsequence converges to a regular solution $u$ of mass $\sigma$ and
\begin{align}
 u_{i,n}&\longrightarrow u_i
 &&\text{in }C^1_{\rm loc}(\R^2)
 \text{ and }C^2_{\rm loc}(\R^2\setminus\{0\}),
 \label{cmp:function-convergence}\\
 f_{i,n}&\longrightarrow\e^{u_i}
 &&\text{in }L^1(\R^2)\text{ and }L^3_{\rm loc}(\R^2).
 \label{cmp:density-convergence}
\end{align}
\begin{equation}\label{cmp:moment-convergence}
 \int_{\R^2}\mathcal K(x)f_{i,n}(x)\dd x
 \longrightarrow\int_{\R^2}\mathcal K(x)\e^{u_i(x)}\dd x
 \qquad(i\in I).
\end{equation}
If \eqref{mom:fixed-nc} holds for $U^\sigma$ and $b$, then $u=U^\sigma$.
\end{theorem}

In particular, the conclusion $u=U^\sigma$ holds when $\sigma$
belongs to a compact set $K\subset\Gamma_0$ and $b$ belongs to a
compact set satisfying \eqref{eq:nc-uniform} on $K$.

The main purpose of this section is to prove
Theorem~\ref{cmp:compactness}.
We first establish several preliminary results.
Lemma~\ref{lem:rectangle} will be used twice: first for the masses
concentrated at a point and then for the masses of the limiting
densities. For $\sigma\in\Gamma_0$, it gives
\[
 0\le\xi\le\sigma,\qquad \xi\ne0,\qquad \Lambda_I^0(\xi)=0
 \quad\Longrightarrow\quad \xi=\sigma.
\]
In the second application, some limiting densities may vanish
identically. The required integral identity must therefore remain
valid when some density components are zero.

\subsection{Local estimates and an integral identity}

We first prove local estimates away from concentration points. We then establish an interaction identity, which will be used before identifying the limiting densities.

Discarding finitely many terms, assume $|\beta_{i,n}|\le1/2$. Fix constants
\begin{equation}\label{cmp:threshold}
 M=\sup_n\sum_i\int_{\R^2}f_{i,n}<\infty,
 \qquad a_* =\max_{i,j}a_{ij}>0,
 \qquad \varepsilon_* =\frac{\pi}{8a_*}.
\end{equation}
Notice that $M<\infty$ follows from $\sigma_n\to\sigma$ and the definition \eqref{eq:mass}. For each fixed $n$, the densities are locally in $L^3$, including at zero, since $3\beta_{i,n}\ge-3/2>-2$. The corresponding local bounds may depend on $n$. The following lemma gives uniform estimates when the mass of a disk is small.

\begin{lemma}\label{cmp:small-mass}
For a fixed disk $D=B_{4r}(p)$ satisfying
\[
 \sum_i\int_D f_{i,n}\,\dd x\le\varepsilon_*,
\]
there exists $C=C(A,N,M,r,p)$, independent of $n$, such that
\begin{equation}\label{cmp:small-mass-estimates}
 \|f_{i,n}\|_{L^3(B_{2r}(p))}\le C,
 \qquad
 \sup_{B_r(p)}u_{i,n}\le C.
\end{equation}
\end{lemma}

\begin{proof}
To prove \eqref{cmp:small-mass-estimates}, we decompose $u_{i,n}$ into the solution of a Dirichlet problem with zero boundary values and a harmonic function. We estimate the two terms separately.

\medskip\noindent\emph{Step 1. Exponential integrability of the potential.}
Let $g_{i,n}=\sum_j a_{ij}f_{j,n}$ and
$\delta_{i,n}=\int_Dg_{i,n}$.  Nonnegativity of $A$ gives
\[
 0\le\delta_{i,n}\le a_*\varepsilon_* =\pi/8.
\]
Let $v_{i,n}$ be the solution of
$-\Delta v_{i,n}=g_{i,n}$ on $D$ with $v_{i,n}|_{\partial D}=0$. By the Green representation and the
maximum principle, we obtain $v_{i,n}\ge0$. The Green function of this
disk satisfies
\[
 0\le G_D(x,y)\le\frac1{2\pi}\log\frac{8r}{|x-y|}.
\]
Therefore, for $s>0$ with $s\delta_{i,n}\le2\pi$ and $\delta_{i,n}>0$,
Jensen's inequality gives
\[
 \e^{sv_{i,n}(x)}
 \le\frac1{\delta_{i,n}}\int_D
 \left(\frac{8r}{|x-y|}\right)^{s\delta_{i,n}/(2\pi)}
 g_{i,n}(y)\,\dd y.
\]
Since $|x-y|\le8r$ and $s\delta_{i,n}/(2\pi)\le1$,
\begin{equation}\label{cmp:green-exponential}
 \int_D\e^{sv_{i,n}}\,\dd x
 \le\sup_{y\in D}\int_D\frac{8r}{|x-y|}\,\dd x
 \le128\pi r^2.
\end{equation}
If $\delta_{i,n}=0$, then $v_{i,n}=0$ and the same bound holds.
In particular, $16\delta_{i,n}\le2\pi$ yields
\begin{equation}\label{cmp:potential-exponential}
 \int_D\e^{16v_{i,n}}\,\dd x\le C(r).
\end{equation}

\medskip\noindent\emph{Step 2. The harmonic remainder.}
The difference $h_{i,n}=u_{i,n}-v_{i,n}$ is harmonic on $D$.
Boundedness of $h_{i,n}$ near the origin makes its harmonic singularity removable for each fixed $n$.  Since $v_{i,n}\ge0$, one has
$h_{i,n}^+\le u_{i,n}^+$.  Almost everywhere,
\[
 u_{i,n}^+
 \le (\log f_{i,n})^+
       +\tfrac12|\log|x||
 \le f_{i,n}+\tfrac12|\log|x||.
\]
Consequently, we get
\[
 \int_D h_{i,n}^+\,\dd x
 \le M+\frac12\int_D|\log|x||\,\dd x.
\]
The last logarithmic integral is finite, including when $0\in D$.
The mean-value estimate for the positive part of a harmonic function
now gives $h_{i,n}\le C$ on $B_{2r}(p)$.  Together with
\eqref{cmp:potential-exponential}, this gives
\[
 \int_{B_{2r}(p)}\e^{16u_{i,n}}\,\dd x\le C.
\]

\medskip\noindent\emph{Step 3. Density and supremum estimates.}
Next, H\"older's inequality with exponents $16/3$ and $16/13$ gives
\begin{align*}
 \int_{B_{2r}(p)}f_{i,n}^3\,\dd x
 &\le
 \left(\int_{B_{2r}(p)}\e^{16u_{i,n}}\,\dd x\right)^{3/16}
 \left(\int_{B_{2r}(p)}|x|^{48\beta_{i,n}/13}\,\dd x\right)^{13/16}.
\end{align*}
Since the disk is fixed and the exponent ranges over a compact interval, its power integral is uniformly bounded.  This proves the
$L^3$ estimate.

For the supremum bound, use the original potential $v_{i,n}$.
For $x\in B_r(p)$, split its Green representation into
$B_{2r}(p)$ and $D\setminus B_{2r}(p)$. H\"older's inequality
on the first part and $|x-y|\ge r$ on the second, give
\[
 v_{i,n}(x)
 \le C(r)\|g_{i,n}\|_{L^3(B_{2r}(p))}
       +\frac{\log8}{2\pi}\int_D g_{i,n}\,\dd x
 \le C.
\]
The harmonic bound from Step~2 now gives $u_{i,n}=v_{i,n}+h_{i,n}\le C$
on $B_r(p)$, with the stated dependence of the constant.
\end{proof}

\begin{lemma}\label{cmp:interaction}
For $X\in C_c^1(\R^2,\R^2)$, set
\[
 D_X(x,y)=
 \frac{(X(x)-X(y))\cdot(x-y)}{|x-y|^2}
 \qquad(x\ne y).
\]
Every finite-mass solution with $\beta_i>-1$ satisfies
\begin{equation}\label{cmp:interaction-eq}
 \sum_i\int_{\R^2}
 \bigl(\operatorname{div}X+\beta_iX\cdot\mathcal K\bigr)f_i
 =\frac1{4\pi}\sum_{i,j}a_{ij}
 \iint_{\R^2\times\R^2}D_X(x,y)f_i(x)f_j(y)\,\dd x\,\dd y.
\end{equation}
The double integral is absolutely convergent with $|D_X(x,y)|\le\operatorname{Lip}(X)$.
\end{lemma}

\begin{proof}
On $\R^2\setminus\{0\}$, we have
\[
 \operatorname{div}(Xf_i)
 =(\operatorname{div}X+\beta_iX\cdot\mathcal K)f_i
   +f_iX\cdot\nabla u_i.
\]
All terms on the right are integrable by Lemma~\ref{mom:origin}.
Choose $R>0$ such that $\supp X\subset B_R$.
The boundary term on $\partial B_R$ vanishes and
\[
 \left|
 \int_{\partial B_\varepsilon}
 f_iX\cdot\frac{x}{\varepsilon}\,\dd S
 \right|
 \le C\varepsilon^{1+\beta_i}\longrightarrow0
 \qquad(\varepsilon\downarrow0).
\]
Integrating over $B_R\setminus\overline{B_\varepsilon}$
and letting $\varepsilon\downarrow0$ gives
\[
 \int_{\R^2}
 (\operatorname{div}X+\beta_iX\cdot\mathcal K)f_i\,\dd x
 =-\int_{\R^2}f_iX\cdot\nabla u_i\,\dd x.
\]

The estimate in the proof of Lemma~\ref{cmp:gradient} gives
\[
\begin{aligned}
 \iint_{\R^2\times\R^2}
 |X(x)|\,&|\mathcal K(x-y)|f_i(x)f_j(y)\,\dd x\,\dd y\\
 &\qquad\le
 \|X\|_{L^\infty(\R^2)}\|f_i\|_{L^1(\R^2)}
 \bigl(C_p\|f_j\|_{L^p(\R^2)}+\|f_j\|_{L^1(\R^2)}\bigr)
 <\infty.
\end{aligned}
\]
Substituting \eqref{cmp:gradient-representation}, interchanging the integrals and
summing over $i$, we obtain
\[
\begin{aligned}
 &\sum_i\int_{\R^2}
 (\operatorname{div}X+\beta_iX\cdot\mathcal K)f_i\,\dd x\\
 &\quad=\frac1{2\pi}\sum_{i,j}a_{ij}
 \iint_{\R^2\times\R^2}
 X(x)\cdot\mathcal K(x-y)f_i(x)f_j(y)\,\dd x\,\dd y\\
 &\quad=\frac1{4\pi}\sum_{i,j}a_{ij}
 \iint_{\R^2\times\R^2}
 [X(x)-X(y)]\cdot\mathcal K(x-y)
 f_i(x)f_j(y)\,\dd x\,\dd y\\
 &\quad=\frac1{4\pi}\sum_{i,j}a_{ij}
 \iint_{\R^2\times\R^2}
 D_X(x,y)f_i(x)f_j(y)\,\dd x\,\dd y.
\end{aligned}
\]
For the second equality, we interchange $(i,x)$ with $(j,y)$
and use $a_{ij}=a_{ji}$ and
$\mathcal K(y-x)=-\mathcal K(x-y)$.
This proves \eqref{cmp:interaction-eq}.
Finally,
\[
 |D_X(x,y)|
 \le\frac{|X(x)-X(y)|}{|x-y|}
 \le\operatorname{Lip}(X)
 \qquad(x\ne y).
\]
\end{proof}
\begin{remark}\label{cmp:total-pohozaev}
Take
$X_R(x)=x\chi(x/R)$, where $\chi\in C_c^\infty(\R^2)$ equals
$1$ on $B_1$ and $0$ outside $B_2$.  Its Lipschitz constant is
bounded independently of $R$,
$\operatorname{div}X_R\to2$, $X_R\cdot\mathcal K\to1$ and
$D_{X_R}(x,y)\to1$ for $x\ne y$.  Dominated convergence in
\eqref{cmp:interaction-eq} gives
\[
 \sum_i(2+\beta_i)\int f_i
 =\frac1{4\pi}\sum_{i,j}a_{ij}
   \left(\int f_i\right)\left(\int f_j\right).
\]
This gives $\Lambda_I^\beta(\sigma)=0$ with the convention
\eqref{eq:mass}. 
\end{remark}

\subsection{Exclusion of concentration}

The next proof uses the interaction identity before identifying the limiting densities. We choose $X_r(x)=x-p$ near the point $p$ under consideration. This choice fixes the value of the pair kernel there and prevents a defect on the diagonal of the limiting product measure.

\begin{proposition}\label{cmp:atoms}
Let $u_n$ be finite-mass solutions of \eqref{eq:system} such that
\[
 \begin{gathered}
 \beta_n\to0,\qquad \sigma_n\to\sigma\in\Gamma_0,\qquad
 f_{i,n}\,\dd x\rightharpoonup\mu_i 
 \end{gathered}
\]
for some finite nonnegative Radon measures $\mu_{i}\in \mathcal{M}_+(\R^2)$. Here, $f_{i,n}=|x|^{\beta_{i,n}}\e^{u_{i,n}}$ for $i\in I$. Then the following assertions hold.
\begin{enumerate}
\renewcommand{\labelenumi}{(\roman{enumi})}

\item For $p\in\R^2$, set $\xi_i=\mu_i(\{p\})/(2\pi)$. Then
\begin{equation}\label{cmp:atomic-pohozaev}
 4\sum_i\xi_i-\sum_{i,j}a_{ij}\xi_i\xi_j=0.
\end{equation}
If $\xi\ne0$, then
\[
 \xi=\sigma,\qquad
 \mu_i=2\pi\sigma_i\delta_p\quad(i\in I).
\]

\item If \eqref{cmp:median} holds, then
\[
 \left\{p\in\R^2:\sum_i\mu_i(\{p\})>0\right\}
 \subset\partial B_1,\qquad
 \mu_i(\{0\})=0\quad(i\in I).
\]
\end{enumerate}
\end{proposition}

\begin{proof}
\emph{Step 1. Decomposition of the limiting measures.}
The convergence of $\sigma_n$ gives a uniform total mass bound.
Let $\mu=\sum_i\mu_i$ and set
\[
 S=\{p\in\R^2:\mu(\{p\})\ge\varepsilon_*\}.
\]
The set $S$ is finite. For $p\notin S$, choose $r>0$ such that
\[
 \mu(B_{4r}(p))<\varepsilon_*,
 \qquad
 \mu(\partial B_{4r}(p))=0.
\]
By weak convergence, Lemma~\ref{cmp:small-mass} applies for all
sufficiently large $n$.
Its uniform $L^3$ estimate and weak compactness in $L^3$ show that
$\mu_i$ has an $L^3$ density on $B_{2r}(p)$.
Hence
\begin{equation}\label{cmp:measure-decomposition}
 \mu_i=\sum_{q\in S}2\pi\xi_{i,q}\delta_q+\rho_i\,\dd x,
 \qquad
 \rho_i\in L^1(\R^2)\cap L^3_{\rm loc}(\R^2\setminus S),
\end{equation}
where $\xi_{i,q}=\mu_i(\{q\})/(2\pi)$.
It remains to prove \eqref{cmp:atomic-pohozaev} for $p\in S$.

\medskip\noindent
\emph{Step 2. Passage to the limit in the integral identity.}
Fix $p\in S$ and choose $r>0$ such that $B_{3r}(p)$ contains no
other point of $S$.
If $p\ne0$, also require $2r<|p|$.
Choose $\chi\in C_c^\infty(\R^2)$ with $0\le\chi\le1$,
$\chi=1$ on $B_1$ and $\chi=0$ outside $B_2$ and set
\[
 X_r(x)=(x-p)\chi((x-p)/r).
\]
The Lipschitz constants are bounded independently of $r$.
We also have
\begin{equation}\label{cmp:localized-kernel-values}
 \begin{aligned}
 D_{X_r}(x,y)&=1
 &&(x,y\in B_r(p),\ x\ne y),\\
 D_{X_r}(x,y)&=0
 &&(x,y\notin B_{2r}(p),\ x\ne y).
 \end{aligned}
\end{equation}
Here, for a vector field $X$, we set
\[
 D_X(x,y)
 :=\frac{(X(x)-X(y))\cdot(x-y)}{|x-y|^2},
 \qquad x,y\in\R^2,\quad x\ne y.
\]

Apply \eqref{cmp:interaction-eq} with this fixed $X_r$.
The term containing $\beta_{i,n}$ tends to zero by the total mass bound.
Indeed, $X_r\cdot\mathcal K=\chi(x/r)$ when $p=0$.
When $p\ne0$, the support of $X_r$ is separated from zero, so
$X_r\cdot\mathcal K$ is again bounded.
The divergence term converges by local weak convergence.

For the double integral, local weak convergence gives
\[
 f_{i,n}(x)f_{j,n}(y)\,\dd x\,\dd y
 \rightharpoonup\mu_i\otimes\mu_j
 \quad\text{locally in }\R^2\times\R^2.
\]
Here $\mu_i\otimes\mu_j$ denotes the product measure, defined by $(\mu_i\otimes\mu_j)(E\times F)=\mu_i(E)\mu_j(F)$ for Borel sets $E,F\subset\R^2$. This follows first for test functions of the form
$\varphi(x)\psi(y)$ and then by uniform approximation on compact sets.
Extend $D_{X_r}$ to the diagonal by
\[
 D_{X_r}(x,x)=\frac12\operatorname{div}X_r(x).
\]
The extension equals $1$ near $(p,p)$ and $0$ near $(q,q)$ for
every $q\in S\setminus\{p\}$.
It is continuous off the diagonal.
Its remaining discontinuities belong to
$\{(x,x):x\notin S\}$, which has zero $\mu_i\otimes\mu_j$ measure
by \eqref{cmp:measure-decomposition}.
Thus the integrals converge after multiplication by a continuous
compactly supported cutoff.

To remove the cutoff, take $L>2(|p|+2r)$.
If $|x|>L$ or $|y|>L$ and $D_{X_r}(x,y)\ne0$, the other variable
belongs to $\supp X_r$, so $|D_{X_r}(x,y)|\le C_r/L$.
Consequently,
\[
 \iint_{\{|x|>L\}\cup\{|y|>L\}}
 |D_{X_r}(x,y)|f_{i,n}(x)f_{j,n}(y)\,\dd x\,\dd y
 \le\frac{C_rM^2}{L}.
\]
The same estimate holds for $\mu_i\otimes\mu_j$.
Passing first to the limit in $n$ and then letting $L\to\infty$,
we obtain
\begin{equation}\label{cmp:limiting-localized-interaction}
 \sum_i\int_{\R^2}\operatorname{div}X_r\,\dd\mu_i
 =\frac1{4\pi}\sum_{i,j}a_{ij}
 \iint_{\R^2\times\R^2}
 D_{X_r}(x,y)\,\dd\mu_i(x)\,\dd\mu_j(y).
\end{equation}

\medskip\noindent
\emph{Step 3. The local mass identity.}
Write
\[
 \mu_i=2\pi\xi_i\delta_p+\nu_i,
 \qquad \nu_i(\{p\})=0.
\]
Since $\operatorname{div}X_r(p)=2$, we have
\[
 \sum_i\int_{\R^2}\operatorname{div}X_r\,\dd\mu_i
 =4\pi\sum_i\xi_i
   +O\left(\sum_i\nu_i(B_{2r}(p))\right).
\]
Moreover, $D_{X_r}(p,p)=1$ and
$D_{X_r}(p,y)=0$ for $y\notin B_{2r}(p)$.
Using \eqref{cmp:localized-kernel-values} and the uniform bound on
$D_{X_r}$, we obtain
\[
\begin{aligned}
 \Bigg|\iint_{\R^2\times\R^2}
 D_{X_r}(x,y)\,&\dd\mu_i(x)\,\dd\mu_j(y)
 -4\pi^2\xi_i\xi_j\Bigg|\\
 &\qquad\le CM\bigl(\nu_i(B_{2r}(p))+\nu_j(B_{2r}(p))\bigr)
 \longrightarrow0
 \qquad(r\downarrow0).
\end{aligned}
\]
Letting $r\downarrow0$ in \eqref{cmp:limiting-localized-interaction}
therefore gives
\[
 4\pi\sum_i\xi_i=\pi\sum_{i,j}a_{ij}\xi_i\xi_j,
\]
which proves \eqref{cmp:atomic-pohozaev}.

\medskip\noindent
\emph{Step 4. Determination of the masses and the point $p$.}
Nonnegativity and convergence of the total component masses give
$0\le\xi_i\le\sigma_i$.
If $\xi\ne0$, Lemma~\ref{lem:rectangle} and
\eqref{cmp:atomic-pohozaev} imply $\xi=\sigma$.
Hence
\[
 2\pi\sigma_i=\mu_i(\{p\})
 \le\mu_i(\R^2)\le2\pi\sigma_i,
\]
so $\mu_i=2\pi\sigma_i\delta_p$ for every $i\in I$.

Suppose now that \eqref{cmp:median} holds.
If $|p|\ne1$, then $\mu_i(\partial B_1)=0$ and weak convergence gives
\[
 \lim_{n\to\infty}\sum_i\int_{B_1}f_{i,n}\,\dd x
 =\begin{cases}
 2\pi\sum_i\sigma_i,&|p|<1,\\
 0,&|p|>1.
 \end{cases}
\]
Both cases contradict \eqref{cmp:median}.
Thus $|p|=1$, which proves the second assertion.
\end{proof}

\begin{lemma}\label{cmp:moment-limit}
Let $h_{i,n}\ge0$ satisfy
\[
 \sup_n\sum_i\int_{\R^2}h_{i,n}(x)\,\dd x<\infty,
 \qquad
 h_{i,n}\,\dd x\rightharpoonup\nu_i,
\]
for some finite nonnegative Radon measures
$\nu_i\in\mathcal M_+(\R^2)$.
Assume also that
\begin{equation}\label{cmp:weighted-ui}
\omega(r):=\sup_n\sum_i\int_{B_r}
 \frac{h_{i,n}(x)}{|x|}\,\dd x\longrightarrow0
 \quad(r\downarrow0).
\end{equation}
Then $\nu_i(\{0\})=0$ and 
\begin{equation}\label{cmp:weak-moment-convergence}
 \int_{\R^2}\mathcal K(x)h_{i,n}(x)\,\dd x
 \longrightarrow
 \int_{\R^2}\mathcal K(x)\,\dd\nu_i(x).
\end{equation}
\end{lemma}

\begin{proof}
Fix $r>0$. Choose a smooth nondecreasing function $\theta:\R\to[0,1]$ such that $\theta(s)=0$ for $s\le1$ and  $\theta(s)=1$ for $s\ge2$.
Set
\[
 \chi_k(x)=\theta(k|x|)\theta\bigl(k(r-|x|)\bigr),
 \qquad k\ge1.
\]
Then $\chi_k\in C^\infty_c(\{0<|x|<r\})$ and
\[
 0\le\chi_k(x)\le\chi_{k+1}(x)\le1,\qquad
 \lim_{k\to+\infty}\chi_k(x)=1
 \qquad(0<|x|<r).
\]
For each fixed $k$, the function $\chi_k(x)/|x|$,
extended by zero at the origin, is continuous and compactly
supported. Local weak convergence therefore gives
\[
\begin{aligned}
 \sum_i\int_{\R^2}\frac{\chi_k(x)}{|x|}\,\dd\nu_i(x)
 &=
 \lim_{n\to\infty}
 \sum_i\int_{\R^2}
 \frac{\chi_k(x)}{|x|}h_{i,n}(x)\,\dd x\\
 &\le
 \sup_n\sum_i\int_{B_r}
 \frac{h_{i,n}(x)}{|x|}\,\dd x
 =\omega(r).
\end{aligned}
\]
Letting $k\to\infty$ and applying the monotone convergence
theorem, we obtain
\[
 \sum_i\int_{0<|x|<r}\frac{\dd\nu_i(x)}{|x|}
 \le\omega(r).
\]
The corresponding mass estimate follows directly from
$\sum_i\int_{B_r}h_{i,n}\le r\omega(r)$, by continuous cutoff
approximation and an arbitrarily small enlargement of the disk.
It implies $\sum_i\nu_i(\{0\})=0$.  The displayed weighted
bound therefore includes the whole disk.

Denote $M_0:=\sup_n\sum_i\int_{\mathbb R^2}h_{i,n}<\infty$.  For every $R>0$,
\[
 \sum_i\int_{|x|>R}\frac{h_{i,n}(x)}{|x|}\,\dd x
 \le\frac{M_0}R,
 \qquad
 \sum_i\int_{|x|>R}\frac{\dd\nu_i(x)}{|x|}
 \le\frac{M_0}R.
\]
Choose a smooth cutoff function $\chi_{r,R}\in C_c^\infty([r/2,2R])$ with $0\le\chi_{r,R}\le1$ and
equals to $1$ on $\{r\le|x|\le R\}$. Since $\chi_{r,R}\mathcal K\in C^\infty_c(\R^2,\R^2)$,
local weak convergence gives
\[
 \int_{\R^2}\chi_{r,R}(x)\mathcal K(x)h_{i,n}(x)\,\dd x
 \longrightarrow
 \int_{\R^2}\chi_{r,R}(x)\mathcal K(x)\,\dd\nu_i(x).
\]
Moreover, $1-\chi_{r,R}$ vanishes on $\{r\le|x|\le R\}$. Since $|\mathcal K(x)|=|x|^{-1}$,
the above estimates imply
\[
\begin{aligned}
 \int_{\R^2}(1-\chi_{r,R})|\mathcal K|h_{i,n}\,\dd x
 &\le\omega(r)+\frac{M_0}{R},\\
 \int_{\R^2}(1-\chi_{r,R})|\mathcal K|\,\dd\nu_i
 &\le\omega(r)+\frac{M_0}{R}.
\end{aligned}
\]
In particular, we get
\begin{align*}
 \left|\int\mathcal K h_{i,n}\,\dd x-\int\mathcal K\,\dd\nu_i\right|\le
 \left|\int\chi_{r,R}\mathcal K h_{i,n}\,\dd x
       -\int\chi_{r,R}\mathcal K\,\dd\nu_i\right|+2\omega(r)+\frac{2M_0}{R}.
\end{align*}
Take $n\to\infty$ at fixed $r,R$, then let
$r\downarrow0$ and $R\to\infty$. This proves
\eqref{cmp:weak-moment-convergence}.
\end{proof}

\begin{corollary}\label{cmp:no-atoms}
Under the assumptions and normalization of Theorem~\ref{cmp:compactness}, suppose
$f_{i,n}\,\dd x\rightharpoonup\mu_i$ locally. Then
\[
 \mu_i(\{0\})=0,\qquad
 \eqref{cmp:weak-moment-convergence}\ \text{holds with }h_{i,n}=f_{i,n},
\]
\[
 b\cdot\sigma\ne0\quad\Longrightarrow\quad
 \mu_i(\{p\})=0\qquad(p\in\R^2,\ i\in I).
\]
\end{corollary}

\begin{proof}
Proposition~\ref{cmp:atoms} gives $\mu_i(\{0\})=0$ for every
$i\in I$. Choose $r_0\in(0,1)$ such that
\[
 \sum_i\mu_i(B_{4r_0})<\varepsilon_*,
 \qquad
 \sum_i\mu_i(\partial B_{4r_0})=0.
\]
By local weak convergence, Lemma~\ref{cmp:small-mass} applies
for all sufficiently large $n$. Since $\beta_n\to0$, we may assume that
\[
 u_{i,n}\le C_0\quad\text{on }B_{r_0},
 \qquad |\beta_{i,n}|\le\frac12
 \qquad(i\in I).
\]
Consequently, for $0<r<r_0$,
\begin{align}
 \sum_i\int_{B_r}\frac{f_{i,n}(x)}{|x|}\,\dd x
 &\le2\pi N\e^{C_0}\int_0^r s^{-1/2}\,\dd s
 =4\pi N\e^{C_0}r^{1/2}.
 \label{cmp:origin-moment-bound}
\end{align}
Thus \eqref{cmp:weighted-ui} holds with $h_{i,n}=f_{i,n}$.
Lemma~\ref{cmp:moment-limit} yields
\eqref{cmp:weak-moment-convergence} with $\nu_i=\mu_i$.

Dividing \eqref{mom:identity-eq} by $\tau_n>0$, we obtain
\[
 \sum_i b_{i,n}
 \int_{\R^2}\mathcal K(x)f_{i,n}(x)\,\dd x=0.
\]
Using the convergence just proved and $b_n\to b$, we may
let $n\to\infty$ to obtain
\begin{equation}\label{cmp:limit-moment-zero}
 \sum_i b_i\int_{\R^2}\mathcal K(x)\,\dd\mu_i(x)=0.
\end{equation}

Suppose now that $b\cdot\sigma\ne0$ and
$\sum_i\mu_i(\{p\})>0$ for some $p\in\R^2$.
Proposition~\ref{cmp:atoms} gives
\[
 \mu_i=2\pi\sigma_i\delta_p\quad(i\in I),
 \qquad |p|=1.
\]
Substituting into \eqref{cmp:limit-moment-zero}, we find
\[
 0=2\pi(b\cdot\sigma)\mathcal K(p).
\]
This contradicts $b\cdot\sigma\ne0$ and $|p|=1$.
\end{proof}

\begin{remark}
Under the assumptions of Corollary~\ref{cmp:no-atoms},
if $b\cdot\sigma\ne0$, then no concentration occurs at any
point of $\R^2$. Loss of mass at infinity will be excluded
in the proof of Theorem~\ref{cmp:compactness}.
\end{remark}

\subsection{Proof of Theorem~\ref{cmp:compactness}}

We next exclude vanishing of individual components in the limit. For this purpose, we first establish a mass identity using the gradient representation for all $N$ components, allowing some limiting densities to vanish identically.

\begin{proposition}\label{cmp:partial}
Let $W_i\in C^\infty(\R^2)$ and $c_i\ge0$ satisfy
\[
 \rho_i=c_i\e^{W_i}\in L^1(\R^2),\qquad
 \nabla W_i(x)=-\frac1{2\pi}\sum_j a_{ij}
       \int_{\R^2}\mathcal K(x-y)\rho_j(y)\,\dd y.
\]
Then, with $\xi_i=(2\pi)^{-1}\int_{\R^2}\rho_i$,
\begin{equation}\label{cmp:partial-pohozaev}
 \Lambda_I^0(\xi)=4\sum_i\xi_i-
                  \sum_{i,j}a_{ij}\xi_i\xi_j=0.
\end{equation}
\end{proposition}

\begin{proof}
Fix $X\in C_c^1(\R^2;\R^2)$. The calculation in the proof of Lemma~\ref{cmp:interaction}
therefore gives
\[
 \sum_i\int_{\R^2}\rho_i\operatorname{div}X\,\dd x
 =
 \frac1{4\pi}\sum_{i,j}a_{ij}
 \iint_{\R^2\times\R^2}
 D_X(x,y)\rho_i(x)\rho_j(y)\,\dd x\,\dd y.
\]

Choose $\chi$ as in Remark~\ref{cmp:total-pohozaev} and set
$X_R(x)=x\chi(x/R)$. The passage to the limit $R\to\infty$
is the same as in that remark. We omit the details and obtain
\[
 2\sum_i\int_{\R^2}\rho_i\,\dd x
 =
 \frac1{4\pi}\sum_{i,j}a_{ij}
 \left(\int_{\R^2}\rho_i\,\dd x\right)
 \left(\int_{\R^2}\rho_j\,\dd x\right).
\]
Substituting $\int_{\R^2}\rho_i\,\dd x=2\pi\xi_i$
and dividing by $\pi$ proves \eqref{cmp:partial-pohozaev}.
\end{proof}

\begin{proof}[Proof of Theorem~\ref{cmp:compactness}]
Assume first that $b\cdot\sigma\ne0$.

\medskip\noindent
\emph{Step 1. Local estimates and convergence of the gradients.}
We first obtain local bounds for the densities and prove
convergence of the gradients. These estimates will be used
to determine the limits of the components in Step 2.

By the total mass bound, after passing to a subsequence,
\[
 f_{i,n}\,\dd x\rightharpoonup\mu_i
 \quad\text{locally in }\R^2
 \qquad(i\in I),
\]
where $\mu_i$ are finite nonnegative Radon measures.
Corollary~\ref{cmp:no-atoms} gives
\[
 \mu_i(\{p\})=0
 \qquad(p\in\R^2,\ i\in I).
\]
Lemma~\ref{cmp:small-mass} and a finite covering argument
therefore give uniform local $L^3$ bounds for $f_{i,n}$
and uniform local upper bounds for $u_{i,n}$.
By \eqref{cmp:measure-decomposition},
\[
 \mu_i=\rho_i\,\dd x,\qquad
 \rho_i\in L^1(\R^2)\cap L^3_{\rm loc}(\R^2).
\]

We claim that
\begin{equation}\label{cmp:gradient-limit}
 \begin{gathered}
 \nabla u_{i,n}\longrightarrow P_i
 \quad\text{in }C^0_{\rm loc}(\R^2,\R^2),\\
 P_i(x)=-\frac1{2\pi}\sum_j a_{ij}
 \int_{\R^2}\mathcal K(x-y)\rho_j(y)\,\dd y.
 \end{gathered}
\end{equation}
By \eqref{cmp:gradient-representation}, it suffices to
prove local uniform convergence of the corresponding
integrals. We estimate separately the regions near $x=y$
and near infinity.

Let $Q\subset\R^2$ be compact. For $R$ sufficiently large,
\[
 \sup_{x\in Q}
 \int_{|y|>R}|\mathcal K(x-y)|f_{j,n}(y)\,\dd y
 \le\frac{CM}{R}.
\]
For $0<\delta<1$, H\"older's inequality gives
\[
 \sup_{x\in Q}
 \int_{|x-y|<\delta}
 \frac{f_{j,n}(y)}{|x-y|}\,\dd y
 \le
 C_Q\left(\int_{|z|<\delta}|z|^{-3/2}\,\dd z\right)^{2/3}
 \le C_Q\delta^{1/3}.
\]
Both estimates hold with $f_{j,n}$ replaced by $\rho_j$.

Using smooth cutoffs near $x=y$ and outside $B_R$, we obtain
a continuous kernel with compact support in $y$.
Its integrals converge by local weak convergence.
Uniform continuity of the kernel and the total mass bound
give equicontinuity in $x$, so the convergence is uniform
on $Q$.
Letting first $n\to\infty$, then $\delta\downarrow0$ and
$R\to\infty$, we obtain \eqref{cmp:gradient-limit}.
In particular, $P_i$ is continuous on $\R^2$.

\medskip\noindent
\emph{Step 2. Convergence of the components.}
We next identify the limiting densities and verify the
hypotheses of Proposition~\ref{cmp:partial}.
Fix $q\ne0$ and write
\[
 f_{i,n}(x)
 =
 \e^{u_{i,n}(q)}|x|^{\beta_{i,n}}
 \exp\bigl(u_{i,n}(x)-u_{i,n}(q)\bigr).
\]
We therefore consider the point values $u_{i,n}(q)$ and
the differences $u_{i,n}-u_{i,n}(q)$ separately.

Since $u_{i,n}(q)$ is bounded above, after passing to a
further subsequence we may assume that
\[
 u_{i,n}(q)\longrightarrow\ell_i
 \in\R\cup\{-\infty\}
 \qquad(i\in I).
\]
Set $c_i=\e^{\ell_i}$ if $\ell_i$ is finite and $c_i=0$
otherwise.

To determine the limit of the differences, we use
\[
 u_{i,n}(x)-u_{i,n}(q)
 =
 \int_0^1
 \nabla u_{i,n}(q+s(x-q))\cdot(x-q)\,\dd s.
\]
By \eqref{cmp:gradient-limit}, the right-hand side converges
uniformly on compact sets to
\begin{equation}\label{cmp:normal-differences}
 W_i(x)=\int_0^1
 P_i(q+s(x-q))\cdot(x-q)\,\dd s.
\end{equation}
Together with the convergence of the gradients, this gives
$W_i\in C^1(\R^2)$, $\nabla W_i=P_i$, $W_i(q)=0$ and
\begin{equation}\label{cmp:difference-convergence}
 u_{i,n}-u_{i,n}(q)\longrightarrow W_i
 \quad\text{in }C^1_{\rm loc}(\R^2).
\end{equation}
Returning to the decomposition of $f_{i,n}$ and using
$\beta_{i,n}\to0$, we obtain
\[
 f_{i,n}(x)\longrightarrow c_i\e^{W_i(x)}
 \qquad(x\ne0).
\]

We next prove local $L^3$ convergence in order to identify
the measure limit and pass to the limit in the equation.
On each compact set, the upper bound for $u_{i,n}$ and
$|\beta_{i,n}|\le1/2$ imply
\[
 f_{i,n}(x)^3\le C\bigl(1+|x|^{-3/2}\bigr).
\]
The right-hand side is locally integrable.
The dominated convergence theorem gives
\begin{equation}\label{cmp:partial-density-limit}
 f_{i,n}\longrightarrow c_i\e^{W_i}
 \quad\text{in }L^3_{\rm loc}(\R^2).
\end{equation}
Hence $\rho_i=c_i\e^{W_i}$.
Each limiting density is therefore either positive or
identically zero. If $c_i=0$, then $u_{i,n}\to-\infty$
uniformly on compact sets.

Passing to the limit in the equation, we obtain
\[
 -\Delta W_i=\sum_j a_{ij}\rho_j
 \quad\text{in the sense of distributions on }\R^2.
\]
Since $\rho_i=c_i\e^{W_i}$, interior elliptic regularity
gives $W_i,\rho_i\in C^\infty(\R^2)$.
Moreover, $\nabla W_i=P_i$ and \eqref{cmp:gradient-limit}
give the gradient representation in
Proposition~\ref{cmp:partial}.
All its hypotheses are therefore satisfied, including when
some $c_i$ are zero.

\medskip\noindent
\emph{Step 3. Equality of the masses.}
We now prove that every limiting component has the
prescribed mass. Set
\[
 \xi_i=\frac1{2\pi}\int_{\R^2}\rho_i\,\dd x.
\]
Nonnegativity and convergence of the total component masses
give $0\le\xi\le\sigma$.
Proposition~\ref{cmp:partial} gives
$\Lambda_I^0(\xi)=0$.
To apply Lemma~\ref{lem:rectangle}, it remains to show
that $\xi\ne0$.

By \eqref{cmp:median} and local $L^1$ convergence,
\begin{equation}\label{cmp:limit-half-mass}
 \sum_i\int_{B_1}\rho_i\,\dd x
 =\pi\sum_i\sigma_i>0.
\end{equation}
Thus $\xi\ne0$.
Lemma~\ref{lem:rectangle} now yields
\begin{equation}\label{cmp:no-lost-mass}
 \xi=\sigma.
\end{equation}
Since every $\sigma_i$ is positive, all $c_i$ are positive
and all $\ell_i$ are finite.
Consequently,
\[
 u_i=\ell_i+W_i
\]
is a regular solution of \eqref{eq:system} with $\beta=0$
and mass $\sigma$.
By \eqref{cmp:difference-convergence},
$u_{i,n}\to u_i$ in $C^1_{\rm loc}(\R^2)$.

The equality of the masses also gives global $L^1$
convergence of the densities.
Indeed, for each $R>0$, local $L^1$ convergence and
\eqref{cmp:no-lost-mass} imply
\[
 \limsup_{n\to\infty}
 \|f_{i,n}-\rho_i\|_{L^1(\R^2)}
 \le
 2\int_{\R^2\setminus B_R}\rho_i\,\dd x.
\]
Letting $R\to\infty$, we obtain convergence in $L^1(\R^2)$.
Together with \eqref{cmp:partial-density-limit}, this proves
\eqref{cmp:density-convergence}.

On compact subsets of $\R^2\setminus\{0\}$, the factors
$|x|^{\beta_{i,n}}$ converge to $1$ in $C^1$.
Hence the right-hand sides of the equations converge
locally in $C^1$.
Interior Schauder estimates give convergence in
$C^2_{\rm loc}(\R^2\setminus\{0\})$, proving
\eqref{cmp:function-convergence}.
Equation \eqref{cmp:moment-convergence} follows from
Corollary~\ref{cmp:no-atoms}.

\medskip\noindent
\emph{Step 4. The center and scaling.}
It remains to identify the limit under
\eqref{mom:fixed-nc}.
By Theorem~\ref{thm:published-classification}, there exist
$a\in\R^2$ and $\lambda>0$ such that
\begin{equation}\label{cmp:translated-limit}
 u_i(x)=U_i^\sigma(\lambda(x-a))+2\log\lambda.
\end{equation}
We shall prove that $a=0$ and $\lambda=1$.

First, \eqref{cmp:limit-moment-zero} gives
\[
 \sum_i b_i\int_{\R^2}\mathcal K(x)\e^{u_i(x)}\,\dd x=0.
\]
To evaluate this integral, observe that the centered
densities $\lambda^2\e^{U_i^\sigma(\lambda x)}$ have mass
$E_i^\sigma(\lambda r)$ in $B_r$.
If $a\ne0$, Proposition~\ref{mom:center} therefore gives
\[
 0=
 \frac{a}{|a|^2}
 \sum_i b_iE_i^\sigma(\lambda|a|).
\]
This contradicts \eqref{mom:fixed-nc}.
Hence $a=0$.

We next use the normalization to determine $\lambda$.
The reference solution satisfies
\[
 \sum_iE_i^\sigma(1)=\pi\sum_i\sigma_i.
\]
By \eqref{cmp:limit-half-mass} and
\eqref{cmp:translated-limit},
\[
 \sum_iE_i^\sigma(\lambda)=\sum_iE_i^\sigma(1).
\]
Since $r\mapsto\sum_iE_i^\sigma(r)$ is strictly increasing,
we have $\lambda=1$.
Thus $u=U^\sigma$.

Finally, if $\sigma\in K$ and $b$ belongs to the compact
set of vectors in \eqref{eq:nc-uniform},
Corollary~\ref{mom:nc-consequences} gives
$b\cdot\sigma\ne0$.
Condition \eqref{mom:fixed-nc} also follows from
\eqref{eq:nc-uniform}.
The preceding argument therefore gives $u=U^\sigma$.
\end{proof}

\begin{remark}\label{cmp:scope}
The nonzero parameter condition $\tau_n>0$ is used when the
inverse-moment identity is divided by the singularity size.
At $\beta_n=0$ the regular translation family supplies no
such restriction.  The condition $b\cdot\sigma\ne0$
excludes point masses and gives convergence to a regular solution.
The condition at all radii is needed to locate its center.
Neither conclusion alone asserts that the approximating
singular solutions are radial. That additional step follows from the bounded linearized equation
in Section~\ref{sec:rigidity}.
\end{remark}

\section{Radial symmetry and nondegeneracy}\label{sec:rigidity}

We first prove a uniform estimate after inversion. We then determine the bounded solutions of the linearized system. Applying this result to the rotation derivative gives radial symmetry. Throughout this section, $K$ and $\mathcal B$ satisfy the hypotheses of Theorem~\ref{thm:signed}.

\subsection{Uniform estimates at infinity}

Recall the assumptions in Theorem \ref{cmp:compactness}. After taking a subsequence and imposing \eqref{cmp:median}, we write
\begin{equation}\label{rig:compactness-input}
 \begin{gathered}
 \beta_n=\tau_n b_n\to0,\qquad b_n\in\mathcal B,\qquad
 \sigma_n\to\eta\in K,\\
 u_n\longrightarrow U^\eta
 \quad\text{in }C^1_{\rm loc}(\R^2)
 \text{ and }C^2_{\rm loc}(\R^2\setminus\{0\}),\\
f_{i,n}:=|x|^{\beta_{i,n}}\e^{u_{i,n}}       \longrightarrow\e^{U_i^\eta}\quad\text{in }L^1(\R^2).
 \end{gathered}
\end{equation}
Theorem \ref{thm:published-asymptotics} and compactness of $K$ give
\[
 D=\min_{\eta\in K,\ i\in I}\bigl((A\eta)_i-2\bigr)>0,
 \qquad d_*=D/2.
\]
Discarding finitely many terms, we have
\begin{equation}\label{rig:uniform-gaps}
 d_{i,n}:=(A\sigma_n)_i-2-\beta_{i,n}\ge d_*,
 \qquad \gamma_{i,n}\ge\frac32.
\end{equation}
The exponents and masses also have uniform upper bounds.

\begin{lemma}\label{rig:kelvin}
For a sequence satisfying \eqref{rig:compactness-input}, define
\[
 \mathcal I(y)=\frac{y}{|y|^2},\qquad
 \widehat u_{i,n}(y)
    =u_{i,n}(\mathcal I(y))-m_{i,n}\log|y|,\qquad m_n=A\sigma_n.
\]
There exist $R\in(0,1)$ and $C<\infty$, independent of $n$, such that
\begin{equation}\label{rig:kelvin-bound}
 \sup_{i,n}\|\widehat u_{i,n}\|_{L^\infty(B_R)}\le C.
\end{equation}
\end{lemma}

\begin{proof}
\medskip\noindent\emph{Step 1. The equation after inversion.}
We first show that the transformed equation holds on the full
disk, so that the Dirichlet representation can be applied.

The transformation law for the Laplacian gives, away from zero,
\begin{equation}\label{rig:kelvin-equation}
 -\Delta\widehat u_{i,n}
 =\sum_j a_{ij}|y|^{\widehat\beta_{j,n}}
                         \e^{\widehat u_{j,n}},
 \qquad
 \widehat\beta_{i,n}
 =m_{i,n}-\beta_{i,n}-4=d_{i,n}-2.
\end{equation}
Set
$\widehat f_{i,n}=|y|^{\widehat\beta_{i,n}}
\e^{\widehat u_{i,n}}$.
Since inversion has Jacobian $|y|^{-4}$, we also have
\begin{equation}\label{rig:kelvin-mass}
 \int_{B_R}\widehat f_{i,n}\,\dd y
 =\int_{|x|>1/R}f_{i,n}\,\dd x.
\end{equation}

For each fixed $n$, Theorem~\ref{thm:published-asymptotics}
shows that $\widehat u_{i,n}$ is bounded near zero.
Since $\widehat\beta_{i,n}>-2$, the right-hand side of
\eqref{rig:kelvin-equation} belongs locally to $L^p$ for
some $p>1$, which may depend on $n$.
Its local Newton potential is therefore bounded.
Subtracting this potential leaves a bounded harmonic function
on a punctured disk, whose singularity at zero is removable.
Thus \eqref{rig:kelvin-equation} holds in the sense of
distributions on the full disk.

\medskip\noindent\emph{Step 2. Reduction to an $L^p$ estimate.}
We now separate the harmonic part from the Dirichlet potential.
The boundary values control the former, while an $L^p$ bound
for the source, with $p>1$, will control the latter.

Let $R\in(0,1)$ be fixed below and write $D_R=B_R$.
On $\partial D_R$, inversion takes values on the fixed circle
of radius $1/R$.
The local convergence in \eqref{rig:compactness-input}
and boundedness of $m_n$ give
\begin{equation}\label{rig:kelvin-boundary}
 \sup_{i,n}
 \|\widehat u_{i,n}\|_{L^\infty(\partial D_R)}<\infty.
\end{equation}
Let $G_{D_R}$ be the Dirichlet Green function of $-\Delta$
on $D_R$ and define
\[
 V_{i,n}(x)
 =\int_{D_R}G_{D_R}(x,y)
       \sum_j a_{ij}\widehat f_{j,n}(y)\,\dd y,
 \qquad
 H_{i,n}=\widehat u_{i,n}-V_{i,n}.
\]
By Step~1, $H_{i,n}$ is harmonic on $D_R$ and equals
$\widehat u_{i,n}$ on $\partial D_R$.
The maximum principle gives
\begin{equation}\label{rig:harmonic-bound}
 \|H_{i,n}\|_{L^\infty(D_R)}\le C.
\end{equation}
Moreover, $V_{i,n}\ge0$ because its source is nonnegative.

For every $p>1$, the logarithmic Green kernel has a bounded
$L^{p'}(D_R)$ norm in $y$, uniformly for $x\in D_R$,
where $p'=p/(p-1)$.
Hence H\"older's inequality gives
\[
 \|V_{i,n}\|_{L^\infty(D_R)}
 \le C_{p,R}\sum_j a_{ij}
             \|\widehat f_{j,n}\|_{L^p(D_R)}.
\]
It therefore remains to choose $p>1$ and $R>0$ such that
these source norms are bounded independently of $n$.

\medskip\noindent\emph{Step 3. The uniform source estimate.}
Since
$\widehat f_{i,n}
=|y|^{\widehat\beta_{i,n}}\e^{\widehat u_{i,n}}$,
we estimate the power weight and the exponential separately.
Choose
\begin{equation}\label{rig:holder-exponents}
 a=\min\left\{\frac12,\frac{d_*}{2}\right\},
 \qquad
 p=1+\frac a4,
 \qquad
 p_1=1+\frac a2,
 \qquad
 s=\frac{p p_1}{p_1-p}.
\end{equation}
Then $1<p<p_1$ and $1/p=1/p_1+1/s$.
By \eqref{rig:uniform-gaps},
\begin{equation}\label{rig:power-margin}
 p_1\widehat\beta_{i,n}
 \ge\left(1+\frac a2\right)(-2+2a)
 =-2+a+a^2>-2.
\end{equation}
Consequently,
\[
 \int_{D_R}|y|^{p_1\widehat\beta_{i,n}}\,\dd y
 =\frac{2\pi R^{2+p_1\widehat\beta_{i,n}}}
        {2+p_1\widehat\beta_{i,n}}
 \le\frac{2\pi}{a+a^2}.
\]

To estimate the exponential, we choose $R$ so that the
Dirichlet source has sufficiently small mass.
The global $L^1$ convergence in \eqref{rig:compactness-input}
implies
\[
 \lim_{R\downarrow0}\sup_n
 \sum_j\int_{|x|>1/R}f_{j,n}(x)\,\dd x=0.
\]
Thus \eqref{rig:kelvin-mass} and nonnegativity of $A$ allow
us to fix $R\in(0,1)$ such that
\begin{equation}\label{rig:small-kelvin-mass}
 M_{i,n}
 :=\int_{D_R}\sum_j a_{ij}\widehat f_{j,n}\,\dd y
 \le\frac{2\pi}{s}
 \qquad(i\in I,\ n\ge1).
\end{equation}
Applying \eqref{cmp:green-exponential} with
$D=D_R$, $4r=R$ and $v_{i,n}=V_{i,n}$ gives
\[
 \int_{D_R}\e^{sV_{i,n}}\,\dd y\le C(R),
\]
because $sM_{i,n}\le2\pi$.
Together with \eqref{rig:harmonic-bound}, this yields
\begin{equation}\label{rig:exponential-integrability}
 \int_{D_R}\e^{s\widehat u_{i,n}}\,\dd y\le C.
\end{equation}
H\"older's inequality with the exponents in
\eqref{rig:holder-exponents} now gives
\begin{equation}\label{rig:kelvin-source-lp}
 \|\widehat f_{i,n}\|_{L^p(D_R)}
 \le
 \bigl\||y|^{\widehat\beta_{i,n}}\bigr\|_{L^{p_1}(D_R)}
 \|\e^{\widehat u_{i,n}}\|_{L^s(D_R)}
 \le C.
\end{equation}
The estimate in Step~2 therefore bounds $V_{i,n}$ uniformly
in $L^\infty(D_R)$.
Since $\widehat u_{i,n}=V_{i,n}+H_{i,n}$,
\eqref{rig:harmonic-bound} proves \eqref{rig:kelvin-bound}.
\end{proof}

In particular, \eqref{rig:kelvin-equation} and
\eqref{rig:kelvin-bound} give the pointwise bound
\begin{equation}\label{rig:kelvin-density-bound}
 0\le\widehat f_{i,n}(y)\le C|y|^{d_*-2}
 \qquad(0<|y|<R).
\end{equation}
Only this bound will be needed for the linearized system.

\subsection{The bounded solutions of the linearized system}

The following identity holds for every bounded solution of the
linearized system. Its proof uses the same symmetry as
Proposition~\ref{mom:identity}.

\begin{lemma}\label{rig:linearized-moment}
Let $u$ satisfy the assumptions of Lemma~\ref{mom:origin} and let
$\phi$ be a bounded solution of \eqref{eq:linearized}. Then
\begin{equation}\label{rig:linearized-moment-identity}
\sum_i\beta_i\int_{\R^2}\mathcal K(x)f_i(x)\phi_i(x)\,\dd x=0.
\end{equation}
\end{lemma}

\begin{proof}
On $\R^2\backslash\{0\}$, we get
\[
 \nabla(f_i\phi_i)
 =\beta_i\mathcal K f_i\phi_i
   +f_i\phi_i\nabla u_i+f_i\nabla\phi_i.
\]
We shall integrate this identity and show that the integrals
of the last two terms cancel after summation over $i$.
For this purpose, we first obtain an integral representation
of $\nabla\phi_i$.

\medskip\noindent\emph{Step 1. The gradient representation.}
Put $g_i=f_i\phi_i$.
The proof of Lemma~\ref{cmp:gradient} gives
$f_i\in L^1(\R^2)\cap L^p(\R^2)$ for some $p>2$.
Since $\phi$ is bounded, the same holds for $g_i$.
The decay in Theorem~\ref{thm:published-asymptotics} also gives
\[
 \int_{\R^2}\log(2+|y|)|g_i(y)|\,\dd y<\infty.
\]
Define
\[
 v_i(x)=-\frac1{2\pi}\sum_j a_{ij}
 \int_{\R^2}\log\frac{|x-y|}{1+|y|}
                  g_j(y)\,\dd y.
\]
On $|x-y|<1$, H\"older's inequality controls the logarithmic
singularity, while the preceding logarithmic moment controls
the term $\log(1+|y|)$.
On $|x-y|\ge1$, we use
\[
 -\log(1+|y|)
 \le\log\frac{|x-y|}{1+|y|}
 \le\log(1+|x|).
\]
Thus $v_i$ is well-defined and
\[
 |v_i(x)|\le C\bigl(1+\log(2+|x|)\bigr).
\]

Since $-\Delta v_i=\sum_j a_{ij}g_j$, the function
$\phi_i-v_i$ is harmonic on $\R^2$.
Its logarithmic growth and the interior gradient estimate
imply that it is constant.
Differentiating the potential as in
Lemma~\ref{cmp:gradient}, we obtain
\[
 \nabla\phi_i(x)
 =-\frac1{2\pi}\sum_j a_{ij}
 \int_{\R^2}\mathcal K(x-y)
                  f_j(y)\phi_j(y)\,\dd y.
\]
The same $L^1$ and $L^p$ estimates show that
$\nabla\phi_i$ is bounded on $\R^2$.
In particular, $f_i\nabla\phi_i$ is integrable.

\medskip\noindent\emph{Step 2. Integration and cancellation.}
By Lemma~\ref{mom:origin}, boundedness of $\phi$ and Step~1,
all three terms in the initial differential identity are
absolutely integrable.
The boundary terms on $\partial B_r$ and $\partial B_R$
are respectively
\[
 O(r^{1+\beta_i}),
 \qquad
 O(R^{-1-d_i}),
 \qquad
 d_i=m_i-2-\beta_i>0.
\]
Since $\beta_i>-1$, both terms tend to zero as
$r\downarrow0$ and $R\to\infty$.
Integration on $B_R\setminus B_r$ therefore gives
\[
 \beta_i\int_{\R^2}\mathcal K(x)f_i(x)\phi_i(x)\,\dd x
 =-\int_{\R^2}
       \bigl(f_i\phi_i\nabla u_i+f_i\nabla\phi_i\bigr)\,\dd x.
\]

To evaluate the sum on the right, substitute the gradient representations from Lemma~\ref{cmp:gradient} and Step 1. The resulting double integrals are absolutely convergent by the estimate in the proof of Proposition~\ref{mom:identity} and boundedness of $\phi$. Hence
\[
\begin{aligned}
 \sum_i\int_{\R^2}
       \bigl(f_i\phi_i\nabla u_i+f_i\nabla\phi_i\bigr)\,\dd x =-\frac1{2\pi}\sum_{i,j}a_{ij}
 \iint_{\R^2\times\R^2}
 \mathcal K(x-y)f_i(x)f_j(y)
       \bigl[\phi_i(x)+\phi_j(y)\bigr]\,\dd x\,\dd y.
\end{aligned}
\]
Interchanging $(i,x)$ with $(j,y)$ leaves the scalar factors
unchanged and reverses the sign of $\mathcal K(x-y)$.
Since $a_{ij}=a_{ji}$, the double sum is therefore zero.
Summing the preceding component identities proves
\eqref{rig:linearized-moment-identity}.
\end{proof}

\begin{proposition}\label{rig:bounded-kernel}
Let $K$ and $\mathcal B$ satisfy the hypotheses of Theorem~\ref{thm:signed}.
There exist $\epsilon,\delta_0>0$ such that every solution $u$ of
\eqref{eq:system}--\eqref{eq:mass} satisfying
\[
 \beta=\tau b,\qquad 0<\tau<\epsilon,\qquad b\in\mathcal B,
 \qquad \dist(\sigma(u),K)<\delta_0
\]
has bounded linearized kernel
\begin{equation}\label{rig:full-kernel}
 \left\{\phi\in L^\infty(\R^2,\R^N):
 -\Delta\phi_i=\sum_j a_{ij}|x|^{\beta_j}\e^{u_j}\phi_j\right\}
 =\Span\{Z\},\qquad Z_i=x\cdot\nabla u_i+\gamma_i.
\end{equation}
\end{proposition}

\begin{proof}
We may restrict to $|\beta_i|<1/2$.
Since $u$ is bounded near zero, $f_i\in L^3_{\rm loc}(\R^2)$.
For every bounded solution $\phi$ of \eqref{eq:linearized}, local
elliptic regularity therefore gives
\[
 u_i,\ \phi_i\in W^{2,3}_{\rm loc}(\R^2)
 \subset C^{1,1/3}_{\rm loc}(\R^2).
\]
Differentiating the dilation family \eqref{eq:dilation} gives
$Z$ as a solution of \eqref{eq:linearized} in the sense of
distributions. Lemma~\ref{mom:origin} and
Theorem~\ref{thm:published-asymptotics} show that $Z$ is bounded and
\[
 Z_i(0)=\gamma_i>0.
\]
Thus $\Span\{Z\}$ is contained in the kernel for each solution under
consideration.

Suppose that the assertion is false. There are solutions $u_n$ with
\[
 \beta_n=\tau_n b_n,\qquad 0<\tau_n\to0,\qquad
 b_n\in\mathcal B,\qquad \dist(\sigma_n,K)\to0
\]
whose bounded linearized kernels contain a solution outside
$\Span\{Z_n\}$. Subtracting a multiple of $Z_n$ and dividing by the
supremum norm, we obtain bounded solutions $\phi_n$ with
\begin{equation}\label{rig:kernel-normalization}
 \phi_{N,n}(0)=0,\qquad
 \|\phi_n\|_\infty:=\max_i\sup_{\R^2}|\phi_{i,n}|=1.
\end{equation}
Dilation preserves these two conditions, so we impose
\eqref{cmp:median} on $u_n$.
After passing to a subsequence, compactness of $K$ and $\mathcal B$ gives
\[
 \sigma_n\to\eta\in K,\qquad b_n\to b\in\mathcal B.
\]
Corollary~\ref{mom:nc-consequences} gives $b\cdot\eta\ne0$.
Theorem~\ref{cmp:compactness} therefore gives
\eqref{rig:compactness-input} with limit $U^\eta$.

\medskip\noindent\emph{Step 1. The regular limiting equation.}
The local $L^3$ bounds for $f_{i,n}$ and
\eqref{rig:kernel-normalization} give uniform local $W^{2,3}$ bounds
for $\phi_n$. A further subsequence converges in
$C^1_{\rm loc}(\R^2)$ to a bounded solution $\phi$ of
\[
 -\Delta\phi_i=\sum_j a_{ij}\e^{U_j^\eta}\phi_j,
 \qquad \phi_N(0)=0.
\]
Here the equation passes to the limit by
\eqref{cmp:density-convergence}.
Theorem~\ref{thm:published-regular-kernel} gives
\[
 \phi_i=a\bigl(x\cdot\nabla U_i^\eta+2\bigr)
       +\xi\cdot\nabla U_i^\eta,
 \qquad a\in\R,\quad\xi\in\R^2.
\]
Since $\nabla U_i^\eta(0)=0$, evaluation of its $N$th component at zero
gives $2a=0$. Thus only the translation terms remain.

\medskip\noindent\emph{Step 2. The integral identity.}
Lemma~\ref{rig:linearized-moment}, divided by $\tau_n$, gives
\[
 \sum_i b_{i,n}\int_{\R^2}\mathcal K(x)f_{i,n}(x)
                         \phi_{i,n}(x)\,\dd x=0.
\]
For $0<r<1<R$, local boundedness of $u_n$, bounded total masses and
\eqref{rig:kernel-normalization} imply
\[
 \int_{B_r}|\mathcal K|f_{i,n}|\phi_{i,n}|\le Cr^{1/2},
 \qquad
 \int_{\R^2\setminus B_R}|\mathcal K|f_{i,n}|\phi_{i,n}|
       \le\frac C R.
\]
On the intervening annulus both functions converge locally.
First passing to the limit in $n$, then letting $r\downarrow0$ and
$R\to\infty$, yields
\begin{equation}\label{rig:limit-linearized-moment}
 \sum_i b_i\int_{\R^2}\mathcal K(x)\e^{U_i^\eta(x)}
                                      \phi_i(x)\,\dd x=0.
\end{equation}
For $k,l\in\{1,2\}$, polar integration gives
\[
\begin{aligned}
 \int_{\R^2}\frac{x_l}{|x|^2}\e^{U_i^\eta(x)}
                    \partial_{x_k}U_i^\eta(x)\,\dd x
 &=\pi\delta_{kl}\int_0^\infty
                    \e^{U_i^\eta(r)}(U_i^\eta)'(r)\,\dd r\\
 &=-\pi\delta_{kl}\e^{U_i^\eta(0)}.
\end{aligned}
\]
The last equality uses $\e^{U_i^\eta(r)}\to0$ as $r\to\infty$.
Substituting $\phi_i=\xi\cdot\nabla U_i^\eta$ into
\eqref{rig:limit-linearized-moment}, we find
\[
 \left(\sum_i b_i\e^{U_i^\eta(0)}\right)\xi=0.
\]
The coefficient is nonzero by \eqref{eq:central-gap}, so $\xi=0$.
Consequently $\phi_n\to0$ locally in $C^1$.

\medskip\noindent\emph{Step 3. The supremum norm.}
We show that this local convergence contradicts
\eqref{rig:kernel-normalization}.
Set $\widehat\phi_{i,n}(y)=\phi_{i,n}(\mathcal I(y))$.
By Lemma~\ref{rig:kelvin} and \eqref{rig:kelvin-equation},
\[
 -\Delta\widehat\phi_{i,n}
   =\sum_j a_{ij}\widehat f_{j,n}\widehat\phi_{j,n},
 \qquad
 |\Delta\widehat\phi_{i,n}(y)|\le C|y|^{d_*-2}
 \quad(0<|y|<R).
\]
The equation holds on the full disk after a continuous extension at
zero. To see this, choose $p>1$ with $p(d_*-2)>-2$ if $d_*<2$.
The right-hand side is bounded in $L^p(B_R)$, so its local Newton
potential is continuous. The bounded harmonic difference on the
punctured disk has a removable singularity.

Fix $0<\rho<R$ and put
\[
 M_n(\rho)=\max_i\sup_{|y|=\rho}|\widehat\phi_{i,n}(y)|,
 \qquad
 H_\rho(y)=\frac C{d_*^2}
                    \bigl(\rho^{d_*}-|y|^{d_*}\bigr).
\]
Since $-\Delta H_\rho=C|y|^{d_*-2}$ distributionally on $B_\rho$,
the scalar maximum principle applied to
$M_n(\rho)+H_\rho\pm\widehat\phi_{i,n}$ gives
\[
 \max_i\sup_{|x|\ge\rho^{-1}}|\phi_{i,n}(x)|
 \le M_n(\rho)+\frac C{d_*^2}\rho^{d_*}.
\]
For fixed $\rho$, local convergence makes $M_n(\rho)\to0$ and
$\max_i\sup_{|x|\le\rho^{-1}}|\phi_{i,n}(x)|\to0$.
Thus \eqref{rig:kernel-normalization} implies
\[
 1\le\frac C{d_*^2}\rho^{d_*}.
\]
Letting $\rho\downarrow0$ gives a contradiction.
\end{proof}

\subsection{Proof of the main theorems}

\begin{proof}[Proof of Theorem~\ref{thm:signed}]
Choose $\epsilon$ and $\delta_0$ from
Proposition~\ref{rig:bounded-kernel}, with $\epsilon<1/2$.
Let $u$ be a solution satisfying the hypotheses of the theorem.
Differentiating its rotations gives another solution of the linearized
system,
\[
 \psi_i(x)=x_1\partial_{x_2}u_i(x)-x_2\partial_{x_1}u_i(x).
\]
Lemma~\ref{mom:origin} gives a continuous extension with $\psi_i(0)=0$.
Theorem~\ref{thm:published-asymptotics} shows that $\psi$ is bounded
at infinity, since its leading radial gradient term cancels.
Therefore \eqref{rig:full-kernel} gives $\psi=aZ$ for some $a\in\R$.
At the origin this reads $0=a\gamma_N$, so $a=0$.
The angular derivative of every $u_i$ vanishes on the punctured plane.
Thus $u$ is radial about zero.

Theorem~\ref{thm:published-classification} now gives
$\sigma(u)\in\Gamma_\beta$.
The bounded solutions are given by \eqref{rig:full-kernel}.
The same classification theorem gives uniqueness up to dilation when
the exponents and component masses are fixed.
\end{proof}

\begin{proof}[Proof of Theorem~\ref{thm:nonnegative}]
Take
\[
 \mathcal B_+=\left\{b\in\R^N:b_i\ge0,\ \sum_i b_i=1\right\}.
\]
For every $b\in\mathcal B_+$, $\eta\in K$ and $r>0$,
\[
 \left|\sum_i b_iE_i^\eta(r)\right|=\sum_i|b_i|E_i^\eta(r).
\]
Thus \eqref{eq:NC} holds with $\kappa=1$. Every nonzero nonnegative vector $\beta$ can be written as $\beta=\tau b$, with $\tau=\sum_i\beta_i>0$ and $b\in\mathcal B_+$. If $0\le\beta_i<\epsilon/N$, then $\tau<\epsilon$, where $\epsilon$ is supplied by Theorem~\ref{thm:signed} for this set.
That theorem gives all the assertions.
\end{proof}

\section*{Acknowledgements}
The authors are listed in alphabetical order of their surnames, and all authors have made equal contributions. The authors acknowledge the use of AI tools. All mathematical arguments and proofs in the final manuscript were checked and written by the authors.

\section*{Conflict of interest}
The authors declare that they have no conflicts of interest.

\noindent
\textsc{Zetao Cheng}\\
Department of Mathematics and Research Institute for Natural Sciences, College of Natural Sciences, Hanyang University, 222 Wangsimni-ro Seongdong-gu, Seoul 04763, Republic of Korea\\
\textit{Email address: chengzt20@hanyang.ac.kr}

\medskip

\noindent
\textsc{Haoyu Li}\\
Department of Mathematics, Faculty of Science and Technology, University of Macau, Taipa, Macau, China\\
\textit{Email address: hyli1994@hotmail.com}

\medskip

\noindent
\textsc{Lei Zhang}\\
Department of Mathematics, University of Florida, Gainesville, FL 32611-8105, USA
\textit{Email address: leizhang@ufl.edu}

\end{document}